\documentclass[11pt]{article}

\usepackage[margin=1in]{geometry}
\usepackage{setspace}
\AtBeginDocument{
    \setlength{\abovedisplayskip}{3pt}
    \setlength{\belowdisplayskip}{3pt}
    \setlength{\abovedisplayshortskip}{0pt}
    \setlength{\belowdisplayshortskip}{2pt}
    \setlength{\jot}{1pt}
}

\usepackage{amsmath,amssymb,amsthm}
\usepackage{graphicx}
\usepackage{booktabs}
\usepackage{multirow}
\usepackage{natbib}
\usepackage{hyperref}
\usepackage{amsmath}
\usepackage{amssymb}
\usepackage{amsthm}
\usepackage{mathtools}
\usepackage{booktabs}
\usepackage{multirow}
\usepackage{graphicx}
\usepackage{float}
\usepackage{enumitem}
\usepackage{array}
\usepackage{longtable}
\usepackage{url}
\usepackage{hyperref}
\usepackage{xcolor}
\usepackage{authblk}
\usepackage{adjustbox}

\usepackage{titlesec}

\usepackage{mathtools}
\DeclarePairedDelimiter{\ceil}{\lceil}{\rceil}

\usepackage{bbm}
\usepackage{amsmath}
\DeclareMathOperator*{\argmin}{arg\,min}

\usepackage{subcaption}

\newtheorem{lemma}{Lemma}
\newtheorem{proposition}{Proposition}
\newtheorem{corollary}{Corollary}
\newtheorem{remark}{Remark}

\titlespacing*{\section}{0pt}{0.9em}{0.35em}
\titlespacing*{\subsection}{0pt}{0.7em}{0.25em}
\titlespacing*{\subsubsection}{0pt}{0.6em}{0.2em}

\title{Multi-Period Donation Allocation for Food Banks under Stochastic Capacity and Supply}

\author[1]{Irem Sengul Orgut}
\author[2]{Julie S. Ivy}
\author[3]{Reha Uzsoy}
\author[3]{James Wilson}

\affil[1]{Quinlan School of Business, Loyola University Chicago, Chicago, IL, USA}

\affil[2]{Industrial and Operations Engineering, University of Michigan, Ann Arbor, MI, USA}

\affil[3]{Industrial and Systems Engineering, North Carolina State University, Raleigh, NC, USA}

\date{}

\setlist{
    leftmargin=*,
    topsep=3pt,     
    itemsep=0pt,    
    parsep=0pt,     
    partopsep=0pt   
}

\begin{document}

\begin{singlespace}
\maketitle
\end{singlespace}

\vspace{-3em}

\begin{abstract}
Each food bank within the national Feeding America network faces the challenge of distributing donated food equitably across the food-insecure population it serves while minimizing food storage and waste over a multi-period planning horizon. Donated food supply is highly stochastic, and distribution to beneficiaries is done by autonomous local agencies whose capacities are uncertain and vary over time. Based on our long-standing partnership with a regional food bank in the southeastern United States, we address the equitable and effective distribution of food donations to charitable agencies over a finite planning horizon under uncertain capacity and supply. We derive structural properties of this problem and use them to develop efficient heuristics. An extensive numerical study evaluates the performance of the heuristics relative to a stochastic dynamic programming benchmark, which suffers from the curse of dimensionality. One heuristic consistently provides near-optimal solutions and can be easily implemented by food banks. Our results also apply to other problems involving the equitable and effective distribution of scarce resources under stochastic capacity and supply.
\end{abstract}


\vspace{0.5em}

\section{Introduction}
 
According to the United Nations' State of Food Security and Nutrition in the World report, about 673 million people worldwide faced hunger in 2024 \citep{UN-Food}. Food insecurity is defined as ``a household-level economic and social condition of limited or uncertain access to adequate food" \citep{USDA_FI} and affects both developed and developing countries. It is estimated that one in seven people in the United States (US) faces food insecurity \citep{FA-FI}. Food insecurity coexists with substantial food waste. It is estimated that up to 40\% of all food is wasted annually in the US \citep{FoodWaste}. Food waste constitutes the most significant component of landfills, negatively affecting the environment, and presents a lost opportunity to feed food-insecure individuals and families.

Feeding America oversees a network of 200 food banks that provide food assistance to around 49 million people annually by collecting food donations from governmental, corporate, and individual donors, and then distributing the food to its member food banks \citep{FA-FI}. Typically serving multiple counties in a single state, each food bank is a large warehouse that at the beginning of each time period (i) receives donated food from Feeding America, the federal government, its state government, and its local sources; (ii)  determines the amount of donated food to ship to charitable agencies, such as food pantries and soup kitchens; and (iii) stores the remaining food in the warehouse during the current period, incurring a holding (food-storage) cost. As the final point of distribution in this network, the county-level charitable agencies distribute as much of their food as possible to food-insecure individuals and families in their local communities; any undistributed food is at risk of being discarded, incurring an overshipment (wasted-food) cost.

We study how a food bank can \emph{equitably} distribute its donated food across its food-insecure population \emph{over a multi-period time horizon}  when its \textit{food supply and operational capacity for food storage and distribution are uncertain and time-varying\/}. The objective is to minimize the food bank's expected costs of food storage and waste incurred over its time horizon. Our research is motivated by the challenges faced by our long-term collaborator, the Food Bank of Central and Eastern North Carolina (FBCENC), an affiliate of the national Feeding America network and one of the largest food banks in the US\@. FBCENC serves a 34-county region and over 900 charitable agencies distributing food to the population in need.  It faces a complex decision environment due to the large disparity between food supply and demand as well as Feeding America's requirement of documenting the (nearly) equitable distribution of donated food. Since food banks typically measure equity at the county level, \textit{perfect equity} in food distribution is defined as follows:
\begin{quote}
\textit{A food bank achieves perfect equity in the allocation of donated food if, in each period, the proportion of total distributed food allocated to a county is equal to the proportion of the food-insecure population residing in that county. In this case, each county is said to receive its ``fair share'' of the distributed food.}
\end{quote}

 

For each county in a given food bank's service region, the county's operational capacity to receive, store, and distribute donated food to its final beneficiaries is the aggregation of the operational capacities of the charitable agencies in the county and thus is highly uncertain and time-varying.  One component of this  uncertainty and variability is due to events at the county level: (i) new agencies join the county's network; (ii) existing agencies in the county are closed permanently; or (iii) existing agencies in the county are temporarily suspended when they fail to pay their required fees or meet reporting requirements. Other uncertainties arise due to events at the charitable agencies in each county:    
 (i) changes in an agency's long and short-term food-storage capacity caused by failure or acquisition of refrigeration or other essential equipment such as transportation required to pick up and distribute food; (ii) changes in an agency's budget, since each agency is required to pay the food bank certain fees in order to receive donated food; and (iii) workforce changes at the agencies, which are highly dependent on volunteers. 

At the beginning of each period, a given food bank receives its supply of donated food; and its food-storage and county-level food-shipment decisions are made without knowledge of the actual county-level capacities. 
Upon delivery of donated food to the counties, the county-level operational capacities are observed, and two types of costs are incurred: (i) a holding (food-storage) cost at the food bank for any unshipped food since food must be distributed as quickly as possible to avoid spoilage; and (ii) overshipment (wasted-food) costs at the counties that have received food in excess of their observed operational capacities or their fair share.  
 
For a given food bank, future donated food supply is highly uncertain and varies over time. As a result, food allocation decisions must be made with a long-term perspective, since the amount distributed in any period depends not only on current conditions but also on anticipated future supply and capacity. For example, a food bank may choose to distribute less in the current period when future supply shortages are expected. Conversely, it may distribute more aggressively in anticipation of a surge in donations, particularly when postponing distribution could overwhelm the system and place pressure on limited storage and transportation capacity in later periods. Because food banks often operate with limited IT infrastructure and analytical resources, decision-support tools must be easy to implement and straightforward to interpret.

\subsection{Contributions of This Article}

Our contributions in this article are three-fold.

First, we derive a closed-form optimal solution to the single-period food-distribution problem when the food bank's operational capacity is uncertain but its food supply is known, and show that it is equivalent to the newsvendor problem. A myopic heuristic based on the single-period problem provides a lower bound on the optimal shipment quantities for the multi-period problem. We also identify conditions under which the myopic heuristic is optimal for the multi-period problem, allowing it to be treated as a series of single-period problems. However, these conditions are not realistic for food-bank operations; hence, the general multi-period problem must also be addressed.

Second, we formulate the general multi-period food-distribution problem faced by a food bank when both its food supply and operational capacity are uncertain and time-varying as a stochastic dynamic program. However, the DP's computational burden and the food banks' lack of IT resources render its direct implementation impractical. We analyze the structural properties of the general multi-period problem and establish upper and lower bounds on the food bank's optimal shipment (food-distribution) quantities. Using these solutions, we introduce two additional heuristics: (i) the look-ahead heuristic, which uses information about the probability distribution of the supply in the next period; and (ii) the bounding heuristic, which borrows a concept that was introduced in the inventory literature by \citet{shang2003newsvendor}. 

Using real-world data from FBCENC, we conduct extensive numerical experiments to evaluate the performance of the heuristics relative to the DP-based benchmark. We also explore the sensitivity of the expected cost with respect to the food bank's operational capacity, cost parameters, and target ending inventory levels. The look-ahead heuristic outperforms the other two heuristics in almost every real-world scenario and provides near-optimal solutions. 
We use our results to develop useful insights for managers of nonprofits and humanitarian organizations seeking equitable and effective distribution of scarce resources in similarly uncertain environments.


Previous related work is reviewed in Section \ref{sec:lit_review}. In Section \ref{sec:model}, we formulate a food bank's multi-period food-distribution problem with food supply and operational capacity that are uncertain and time-varying; and define all notation, starting with the single-period problem in Section \ref{sec:singleperiod} and extending to the general multi-period problem in Section \ref{sec:multiperiod}. Section \ref{sec:methods} studies the structural properties of the multi-period problem, presents easy-to-implement heuristics, and derives bounds on the value of the optimal solution. Section \ref{sec:experiments} presents computational experiments based on real-world data.  In Section \ref{sec:conclusions} we summarize the implications of this work, including findings, limitations, and future work.

\section{Literature Review} \label{sec:lit_review}


Research in humanitarian operations has mainly concentrated on the topics of disaster relief (short-term), development (long-term), and community empowerment (sustainability) \citep{zobel2016humanitarian}. Disaster relief problems require urgent responses, while long-term development has an ongoing structure. However, both involve trade-offs between equity, efficiency, and/or effectiveness in these environments \citep{savas1978equity, johnson2007community, beamon2008performance}. Examples of studies considering these objectives are \citet{balcik2008last, hwang1999food, altay2013capability, fiorucci2004real, kondaveti2009decision, psaraftis1985tactical, leclerc2012modeling, mandell1991modelling, vitoriano2011multi} and \citet{ozbay2007stochastic}. In this study, we focus on the objectives of equity and effectiveness. We define \emph{perfect equity} as the case where the fraction of the total available food received by each county in the service region is equal to the fraction of the total population in need residing in that county. \citet{marsh1994equity} provide an extensive review of different equity measures. \emph{Effectiveness} is defined as distributing the maximum food supply to the counties as quickly as possible. These two objectives directly conflict with each other \citep{beamon2008performance, orgut2016modeling, karsu2015inequity, berenguer2016modeling}.



Our work is directly related to the expanding literature on food bank operations. Research in this field can be categorized into strategic, tactical, and operational issues. Strategic issues involve large-scale,  long-term decisions that shape the strategy of the humanitarian organization. Examples include determining promotional strategies for food banks to increase donations \citep{ahire2018harvest}, optimally allocating capacities to the agencies in a food bank's network \citep{orgut2016modeling}, matching the counties or agencies in a food bank's service region to food bank warehouses \citep{martins2019redesigning,islam2021modeling}, and strategically allocating additional capacity to food bank warehouses, agencies, or mobile distribution vehicles to improve access and minimize cost \citep{stauffer2022achieving}. Operational-level studies include last-mile distribution problems for scheduling and routing decisions for distributing food, e.g., \citet{bartholdi1983minimal, solak2014stop, davis2014scheduling}, and \citet{eisenhandler2022multi}. Food rescue operations involve collecting and distributing food donations without warehouse storage and can also be considered operational-level problems \citep{lien2014sequential, balcik2014multi, eisenhandler2019humanitarian, eisenhandler2019segment, nair2018scheduling, blackmon2021rapid,alkaabneh2023routing}. Several papers review contributions in this area from varying perspectives \citep{rivera2023systematic,mahmoudi2022decision,sucharitha2021comprehensive,esmaeilidouki2023food}. We will focus our discussion on tactical decision problems as our study centers on the tactical-level problem of the allocation of food from the food bank to the agencies (aggregated into counties) over time. 

Tactical-level problems involve managing inventory levels and allocating resources considering multiple, usually conflicting objectives. \citet{hasnainelicitation} formulate a multi-objective mathematical model to optimally allocate food donations between food bank warehouses and to agencies, considering equity, effectiveness, and efficiency objectives. The authors use a non-interactive preference elicitation algorithm to determine the food bank's preferences (weights) for these three objectives. \citet{gomez2021food} study the allocation of a set of food items by food banks to beneficiaries considering their vulnerability level, nutritional value, and compatibility among received products to maximize the vulnerability-weighted number of beneficiaries served. \citet{buisman2019donation} focus on charitable agencies and coordinate meal and donation planning to optimize the nutritional value of meals served to beneficiaries. \citet{sengul2023equitable} consider the multi-period, multi-product problem of equitably allocating perishable food donations among the counties in their service region while minimizing food waste. None of these studies incorporate the stochastic nature of the food bank resource allocation problems.

Research on food bank operations has also considered stochastic aspects of the problem. \citet{mohan2013improving} examine efficiency and effectiveness in nonprofit food distribution supply chains using a discrete-event simulation model. \citet{lee2017combining} consider collecting excess crops from farmers and distributing them to the food-insecure population. They develop a Markov Decision Process to optimize collection scheduling and allocation to the food-insecure populations, in the face of uncertain supply and labor (volunteers). \citet{alkaabneh2020unified} develop a stochastic approximate dynamic programming model for food banks to allocate food donations to charitable agencies to maximize the total nutritional value under uncertain supply. \citet{fianu2018markov} develop a discrete-time, discrete-space Markov Decision Process model to determine optimal allocation policies of food donations by food banks to counties. The authors evaluate three allocation rules from the production planning literature and aim to pick the best allocation policy under certain conditions. In contrast, we optimize the allocation amounts without limiting our results to specific allocation policies. Our study also differs from these previous studies as we consider uncertain capacities as well as supply and develop easily-implementable policies to be used by food banks for making allocation decisions.

In this paper, we build upon our previous work, \cite{orgut2016modeling}, in which all parameters were assumed to be deterministic. This work introduced the term “bottleneck county” and identified the critical role of the counties' capacity-to-demand ratios in the optimal solution. \citet{sengul2017modeling} incorporate stochastic capacities in a single-period model, maximizing the benefit achieved through food shipments. In practice, the benefit of shipping one pound of food is difficult to estimate. Hence, our model uses costs instead of benefits to provide a more intuitive representation of the problem, with the objective function of cost minimization. In contrast to \citet{sengul2017modeling}, who treat the problem as a two-stage stochastic model with recourse actions after uncertainty is realized, in our model, decisions must be taken before uncertain outcomes are observed and costs are incurred as the random variables are realized. Further, the model and heuristics proposed in our paper consider multiple periods, which reflects the actual operations of the food bank more realistically but introduces significant additional complexity. We also consider the uncertainty in supply, i.e., food donations to the food bank. Since it is not possible to obtain closed-form solutions for the multi-period problem, we formulate a dynamic program to use as benchmark and derive structural properties that we use to develop upper and lower bounds on the optimal shipment amounts as well as several easy-to-implement heuristics.

We also briefly review the literature on multi-echelon stochastic inventory models, as our model structure parallels those problems. For multi-echelon serial stochastic inventory systems \citep{snyder2019fundamentals}, it is known that an echelon base-stock policy is optimal for each stage, assuming no fixed ordering costs for both finite- \citep{clark1960optimal} and infinite-horizon problems \citep{federgruen1986allocation}. \citet{clark1960optimal} show that the optimal base stock levels for each stage can be obtained exactly by solving a set of recursive equations. However, solving these equations requires substantial computational effort. 
Hence, there is a need for easily-implementable heuristics that can provide near-optimal solutions within reasonable computation times. \citet{shang2011simple} reviews heuristics for periodic and continuous review inventory planning problems.

In his pioneering paper, \citet{veinott1965optimal} considers a multi-period, multi-product problem with linear backorder and holding costs under periodic review and the objective of minimizing the expected discounted cost over the infinite horizon. He examines the conditions under which the myopic policy that optimizes each period’s ordering decisions individually will yield the optimal solution for the entire system. \citet{shang2003newsvendor} consider a multi-echelon, single-product inventory system and develop a simple heuristic that uses the structure of the newsvendor solution to compute upper and lower bounds on the optimal order quantities in each period, which it then uses to obtain approximately optimal base stock levels. We use some of these ideas to develop heuristics for the multi-period food shipment problem.

\citet{levi2007approximation} consider a finite horizon, periodic-review, capacitated inventory system with nonstationary and correlated demands evolving over time. They develop an ordering policy that equalizes the expected marginal holding and backorder costs in a period and provide bounds on its worst-case performance. \citet{levi2008approximation} extend this approach to capacitated systems, noting that myopic policies perform worse for capacitated systems than for uncapacitated systems. They develop bounds on the optimal base stock levels by considering the effect of one period’s backlogging cost over the entire horizon. We also develop a myopic heuristic for the optimal food shipment problem over multiple periods and improve its performance by developing bounds on the optimal shipments in each period.

This study contributes to the field of nonprofit and humanitarian operations research by: ($i$) considering the objectives of equity and effectiveness under stochastic capacities and supply over multiple periods, ($ii$) using ideas from the inventory literature to derive structural properties of this problem, ($iii$) developing easily-implementable and interpretable heuristics, and ($iv$) conducting extensive experiments based on real-world data from our collaborating food bank to demonstrate their promising performance.

  \section{Model Formulation}\label{sec:model}
\subsection{Single-Period Problem}\label{sec:singleperiod}
In this section, we formulate the single-period, stochastic food distribution problem and show that it exhibits the structure of the newsvendor problem \citep{arrow1951optimal, khouja1999single, qin2011newsvendor}. In this problem, the goal is to determine the optimal amount of food to be shipped from the food bank to the counties in the food bank's service region in a given period under uncertain county capacities. Supply and demand are assumed to be deterministic and known. We will introduce supply uncertainty for future periods in the multi-period problem presented in Section \ref{sec:multiperiod}. It is typical to assume that demand is known in this setting since demand is usually based on the size of the food-insecure population, which does not usually fluctuate significantly in the short term. In Section \ref{sec:experiments}, we discuss how we estimate demand, following the approach taken by Feeding America and other studies in food bank operations (e.g., \citet{hasnainelicitation,fianu2018markov,mohan2013improving,sengul2017modeling,orgut2016modeling,orgut2018robust}). In this setting, the food bank will make shipments to the counties prior to the county capacities being observed. After food shipments are made, county capacities are observed. The total cost is calculated as the sum of the holding cost for unshipped food at the food bank and the cost of food overshipment to the counties based on the realized capacity and equity requirements. We start by defining the notation that will be used throughout the paper.

\textbf{\emph{Notation}}

\begin{tabular}{ll}
\multicolumn{2}{l}{\textbf{Set and Indices:}}\\
$\mathcal{J}$: & set of counties, $\mathcal{J}=\{j:1,2,...n\}$\\
\end{tabular}

\begin{tabular}{ll}
\multicolumn{2}{l}{\textbf{Parameters:}}\\
$S$: & total food supply (pounds) to the food bank\\  
$D_j$: & demand (pounds) of county $j$, $j \in \mathcal{J}$, where $\Delta=\sum_{j\in \mathcal{J}}D_{j}$ \\
$C_j$: & \textit{stochastic} capacity (pounds) of county $j$, $j \in \mathcal{J}$ 
\end{tabular}

\begin{tabular}{ll}
\multicolumn{2}{l}{\textbf{Cost Coefficients:}}\\
$h$: & holding cost per pound of food remaining at the food bank\\
$c_o$: & overshipment cost per pound for excess food shipped to the counties
\end{tabular}

\begin{tabular}{ll}
\multicolumn{2}{l}{\textbf{Decision Variables:}}\\
$X_j$: & shipment amount (pounds) to county $j$, $j \in \mathcal{J}$\\
$Y_j$: & overshipped amount (pounds) to county $j$, $j \in \mathcal{J}$\\
\end{tabular}
\medskip

\citet{orgut2016modeling} have shown the importance of each county's capacity-to-demand (CD) ratio to the optimal solution for the single-period, deterministic version of the proposed model. We follow them in defining
\begin{center} $\frac{C_{j}}{D_{j}}$: Capacity-to-demand (CD) ratio of county $j$, $j \in \mathcal{J}$  \end{center}

\begin{center} $R = \displaystyle\min_{j \in \mathcal{J}}{\frac{C_{j}}{D_{j}}}$: Minimum CD ratio over all counties in the service region \end{center}

\begin{center} $\mathcal{B}  = \{l \in \mathcal{J}:\frac{C_{l}}{D_{l}}=R$: Set of counties with the minimum CD ratio, i.e. \textit{bottleneck counties}\}\end{center}

A county’s CD ratio represents its ability to serve the food-insecure population residing in that county. The minimum CD ratio is thus the CD ratio of the county with the most restricted ability to serve its population in need, i.e. with the lowest receiving capacity compared to its demand for donated food. \citet{orgut2016modeling} has shown that the optimal solution to the deterministic counterpart of this problem is to set $X_j^{det} = \min\{R D_j,\frac{SD_j}{\Delta}\}$, since $\Delta=\sum_{j\in \mathcal{J}}D_{j}$. Although in practice there is usually a unique bottleneck county due to the capacity and demand values being large numerical values, for the sake of generality, we define set $\mathcal{B}$ such that it can contain multiple counties. In our context, $R$ is a random variable since the capacities are assumed to be stochastic. Throughout this paper, we place a hat “\textasciicircum” over the symbol for a random variable to indicate its observed value; and a tilde “\textasciitilde” over the symbol for any variable to indicate the aggregation of that quantity over all counties in the service region. We also represent vectors in bold such that $\textbf{X} = [X_1,\dots,X_j,\dots,X_n]$ and $\textbf{C} = [C_1,\dots,C_j,\dots,C_n]$. 

The sequence of events in the single period problem is as follows:

\begin{enumerate}
\item A donated supply of $S$ pounds of food is received. 
\item A total of 
\begin{equation}\label{con:supplySPM}
    \tilde{X}={\sum_{j\in \mathcal{J}}X_{j}} \leq S
\end{equation}
pounds of food is shipped to the counties such that each county $j \in \mathcal{J}$ receives a perfectly equitable assignment of 
\begin{equation}\label{con:equitySPM}
    X_{j} = \tilde{X}\frac{D_{j}}{\sum_{j\in \mathcal{J}}D_{j}} \geq 0 \quad \forall j \in \mathcal{J}
\end{equation}
pounds of food.
\item The capacities $\hat{C_j}$ of the counties are observed.
\item The amount of overshipped food is determined based on the observed capacities of the counties, while preserving perfect equity such that
\begin{equation}\label{con:capSPM}
    X_j - \hat{C_j} \leq Y_j \quad \forall j \in \mathcal{J}
\end{equation}
and 
\begin{equation}\label{con:equitySPM2}
    \frac{X_{j}-Y_{j}}{\tilde{X}-\tilde{Y}} - \frac{D_{j}}{\Delta} = 0 \quad \forall j \in \mathcal{J}.
\end{equation}
\item The holding and overshipment costs are computed for the given values of $\mathbf{X}$ and $\mathbf{\hat{C}}$. 
\end{enumerate}

\noindent Although our models work for any defined period, we typically consider a single period to be a week as food banks usually make one shipment to the agencies in their service region per week. 

The objective function minimizes the sum of the holding cost of the leftover food at the food bank and the expected cost of overshipment and can be  formulated as follows: 
\begin{equation}\label{eq:ObjSPMFS}
\min  \quad h(S-\tilde{X}) + \mathbb{E}_\textbf{C}[Q(\textbf{X})]
\end{equation}
where the expectation is calculated based on the overshipment quantities and is represented as:
\begin{equation}\label{eq:ObjSPMSS}
Q(\textbf{X}) = c_o\tilde{Y}.
\end{equation}

We note that if no county's shipment $X_j$ exceeds its observed capacity $\hat{C_j}$, no overshipment occurs, i.e. we get $Y_j =0$ for all $j$, and the overshipment cost will be zero. Overshipment will be observed because (i) some counties, specifically the bottleneck counties, as we will see, have received more food than their capacity, i.e., $X_j > \hat{C_j}$ for $j \in \mathcal{B}$ resulting in an overshipment at those counties; and, (ii) since those bottleneck counties could not receive the food they were allocated due to limited capacity, i.e. $\hat{C_j} < X_j = \tilde{X}\frac{D_{j}}{\sum_{j\in \mathcal{J}}D_{j}}$ for $j \in \mathcal{B}$, the total food that can be distributed with perfect equity decreases, causing overshipments at other counties. 

Our modeling assumptions are as follows:  
\begin{itemize}
    \item The probability density function (p.d.f.) and the cumulative distribution function (c.d.f.) of the minimum CD ratio, $R$, are known and given by $f_R(r)$ and $F_R(r)$ respectively. We discuss how these functions can be estimated in practice in Section \ref{sec:experiments}.
    \item The support of the random variable $R$ is known and given by $[r_L,r_U]$ where $0 < r_L < \infty$ and $r_L < r_U \leq \infty$.
    \item Shipments that have been sent to the counties cannot be reallocated, or sent back to the food bank, reflecting the actual operations of our partner food bank.
    \item Only one shipment decision can be made for each county in the single-period problem.
\end{itemize}



Lemma \ref{lem:recourse} proves that this stochastic model has complete recourse and the overshipment values can be calculated easily based on the shipment values $\tilde{X}$ and the observed minimum CD ratios $\hat{R}$.

\begin{lemma} \label{lem:recourse}
After the minimum CD ratios $\hat{R}$ are observed, the overshipment amounts $\tilde{Y}$ and $Y_{j}$ can be calculated as:
\begin{equation} \label{eq:Ytilde}
    \tilde{Y} = \max\{0,\tilde{X}-\hat{R}\Delta\},
\end{equation}
\begin{equation} \label{eq:Yjt}
    Y_{j} = \max\{0, X_{j} - \hat{R} D_j\}.
\end{equation}

\end{lemma}

The proof is given in Appendix \ref{Appendix:1}.

Proposition \ref{prop:SPMX} provides the optimal structure for the shipment quantities prior to observing capacities of the counties. We provide the proof of Proposition \ref{prop:SPMX} here as we will use these ideas to identify structural properties for the multi-period problem.  

\begin{proposition}
\label{prop:SPMX}
If (i) the p.d.f. $f_R(r)$ is positive and continuous at every $r \in (r_L, r_U)$; and (ii) the cost coefficients $c_o$ and $h$ are positive, then the unique optimal solution to the single-period problem is:
\begin{equation}
\label{eq:Xstar}
    \tilde{X^*} = \min\left\{S,\Delta F_R^{-1}\left(\frac{h}{c_o}\right)\right\}
\end{equation}
Overshipment occurs if and only if all bottleneck counties receive more food than their observed capacities; in this case, all nonbottleneck counties also receive overshipments.
\end{proposition}

\begin{proof} We start by showing that overshipment is observed if and only if bottleneck counties receive more food than their capacities. \eqref{eq:Ytilde} implies that no overshipment will occur, i.e. $\tilde{Y^*}=0$,  if and only if $\tilde{X} \leq R\Delta$. Hence, overshipment is observed if and only if the bottleneck counties all receive more food than their observed capacities. In other words, by \eqref{eq:Ytilde}, $\tilde{Y^*}>0$ if and only if $\tilde{X}>\hat{R}\Delta$, which, using (\ref{con:equitySPM}), can be written as:
\begin{equation}
    \tilde{X} = \frac{X_j\Delta}{D_j}>\hat{R}\Delta.
\end{equation}
For bottleneck counties $b \in \mathcal{B}$, this indicates
\begin{equation}\label{con:BneckProof-mansc}
    \frac{X_b}{D_b}>\hat{R}=\frac{\hat{C}_b}{D_b}\quad \quad \forall b \in \mathcal{B}
\end{equation}
\begin{equation}\label{con:Bneck-mansc}
    X_b>\hat{C}_b.\quad \quad \forall b \in \mathcal{B}
\end{equation}

We now prove that in this situation all the nonbottleneck counties will also receive overshipments. Since $X_b = \frac{D_b \tilde{X}}{\Delta}$, from (\ref{con:Bneck-mansc}), we have
\begin{equation}\label{con:Bneck2-mansc}
   \frac{D_b \tilde{X}}{\Delta} > \hat{C_b} \quad \quad \forall b \in \mathcal{B} \\
\end{equation}
\begin{equation}\label{con:Bneck3-mansc}
   \frac{X_j}{D_j} = \frac{\tilde{X}}{\Delta} > \frac{\hat{C_b}}{D_b}=R \quad \quad \forall b \in \mathcal{B}, \forall j \in \mathcal{J\setminus{B}},   \\
\end{equation}
yielding $X_j > RD_j$ for all $j \in \mathcal{J\setminus{B}}$. Hence $X_j - RD_j$ represents the amount of overshipment to county $j$.
We can thus write the overall objective function as
\begin{equation}\label{con:Z-mansc}
   Z = h(S-\tilde{X}) + \mathbb{E_\textbf{C}}[c_o\tilde{Y}] 
   = h(S-\tilde{X}) + c_o \tilde{X}
    -   c_o \int_{r_L}^{\frac{\tilde{X}}{\Delta}} r \Delta f_R(r) dr -
   c_o \int_{\frac{\tilde{X}}{\Delta}}^{r_U} \tilde{X} f_R(r) dr.
\end{equation}

To incorporate supply constraint (\ref{con:supplySPM}), we construct the Lagrangian function for the problem:
\begin{equation}
    \label{eq:Prop1_proof2}
    L(\tilde{X},\lambda)=h(S-\tilde{X}) + \mathbb{E}_\textbf{C}[c_o\tilde{Y}]-\lambda(\tilde{X}-S)
\end{equation}
for $\tilde{X} \geq 0$ and $\lambda \leq 0$, where $\lambda$ is the Lagrange multiplier (\citet{bazaraa2004nonlinear}, pp. 210–211).  To proceed, we (a) take the partial derivatives of $L(\tilde{X},\lambda)$ with respect to $\tilde{X}$ and $\lambda$ by invoking a generalized version of Leibniz’s rule (\citet{charles1985protter}, chap. 8, Theorem 3); (b) formulate the resulting Karush-Kuhn-Tucker (KKT) second-order sufficient conditions (\citet{bazaraa2004nonlinear}, Theorem 4.4.2); (c) verify the linear-independence constraint qualification at every $\tilde{X} \geq 0$ (\citet{bazaraa2004nonlinear}, p. 244); (d) solve the KKT conditions for the optimal values $\tilde{X^*}$ and $\lambda^*$; and (e) verify that problem is a convex program with a strictly convex objective function having the unique optimal solution $\tilde{X^*}$ found in (d). Thus we have
\begin{align} \label{eq:Prop1_Lag1}
    \frac{\partial L(\tilde{X},\lambda)}{\partial \tilde{X}} = -h +c_o -c_o(1-F_R(\frac{\tilde{X}}{\Delta}))-\lambda = -h +c_o F_R(\frac{\tilde{X}}{\Delta})-\lambda=0,
\end{align}
and
\begin{equation} \label{eq:Prop1_Lag2}
    \frac{\partial L(\tilde{X},\lambda)}{\partial \lambda} = -\tilde{X} + S = 0,
\end{equation}
since by definition $f_R(r)=0$ for $r<r_L$ or $r>r_U$.

Because by assumption $f_R(r)$ is positive on $(r_L,r_U)$ and $c_o$ and $h$ are positive, from (\ref{eq:Prop1_Lag1}) we see that (i) $h + \lambda= c_o F_R(\frac{\tilde{X}}{\Delta})>0$ for each  $\tilde{X} \in (r_L\Delta,r_U\Delta)$; and (ii) $F_R(r)$ is strictly increasing on $(r_L,r_U)$. Thus we obtain the following as the unique solution of (\ref{eq:Prop1_Lag1}): 
\begin{equation} \label{eq:Prop1_Lag3}
    \tilde{X^*} = \Delta F_R^{-1}\left(\frac{h+\lambda}{c_o}\right).
\end{equation}

In this situation, the optimal value of the Lagrange multiplier $\lambda^*$ is equal to the derivative of the objective function  with respect to $\tilde{X}$  evaluated at the optimum $\tilde{X}=\tilde{X^*}$  (\citet{hadley1964nonlinear}, sec. 3-7). Hence if $\lambda < 0$, then by a standard alternative representation of the KKT necessary conditions (\citet{bazaraa2004nonlinear}, p. 166), we must have 
\begin{equation} \label{eq:LagrangeDerLambdaSPM2}
    \tilde{X^*} = S
\end{equation}
which by (\ref{eq:Prop1_Lag1}) leads to the result
\begin{equation} \label{eq:LagrangeDerXSPM3}
    \lambda^* = -h + c_o F_R\left(\frac{S}{\Delta}\right) <  0.
\end{equation}
Hence if $\lambda^*<0$, then by equation (\ref{eq:LagrangeDerXSPM3}) we have 
\begin{equation} \label{eq:LagrangeDerXSPM4}
    S < \Delta F_R^{-1}\left(\frac{h}{c_o}\right).
\end{equation}
On the other hand, if $\lambda^*=0$, then (\ref{eq:Prop1_Lag1}) and (\ref{eq:LagrangeDerXSPM4}) imply that 
\begin{equation} \label{eq:LagrangeDerXSPM5}
    \tilde{X^*} = \Delta F_R^{-1}\left(\frac{h}{c_o}\right) \leq S.
\end{equation}

Combining (\ref{eq:LagrangeDerLambdaSPM2}), (\ref{eq:LagrangeDerXSPM4}) and (\ref{eq:LagrangeDerXSPM5}), we obtain (\ref{eq:Xstar}).

To show uniqueness of the optimal solution, we compute the second derivative of the Lagrangian function (\ref{eq:Prop1_proof2}) with respect to $\tilde{X}$ in order to show that the objective function is a strictly convex function of $\tilde{X}$.  We have
 \begin{equation} \label{eq:Prop1_2Der}
    \frac{\partial^2 L(\tilde{X},\lambda)}{\partial \tilde{X}^2} = \frac{c_o}{\Delta} f_R\left(\frac{\tilde{X}}{\Delta}\right)
\end{equation}

for each $\tilde{X} \in (r_L \Delta,r_U \Delta)$ and $\lambda \leq 0$.  By assumption $c_o$ is positive and $f_R(\frac{\tilde{X}}{\Delta})>0$  for each $\tilde{X} \in (r_L \Delta,r_U \Delta)$; hence the right-hand side of Equation (\ref{eq:Prop1_2Der}) is also positive for each $\tilde{X} \in (r_L \Delta,r_U \Delta)$ and $\lambda \leq 0$. Thus (\ref{eq:Prop1_proof2}) is strictly convex (\citet{bazaraa2004nonlinear}, p. 569); and it follows that equation (\ref{eq:Xstar}) specifies the unique optimal solution for our problem. The proof is completed. 
\end{proof}

We make the following observations based on Proposition \ref{prop:SPMX}.

\begin{remark}\label{obs:1}
Overshipment cost will be zero if and only if $\tilde{X} \leq R\Delta$, which is equivalent to $X_j \leq \hat{C_j}$ for $j \in \mathcal{B}$ i.e. the total shipment to any bottleneck county does not exceed its observed capacity. Overshipment is only observed if the observed bottleneck counties receive more food than their observed capacities. 
\end{remark}

\begin{remark}\label{obs:2}
For (\ref{eq:Xstar}) to be feasible, the holding cost $h$ cannot exceed the cost of overshipment $c_o$, i.e. we must have $h \leq c_o$.
\end{remark}

To explain the intuition behind Remark \ref{obs:2}, assume that the cost to hold one pound of food at the food bank is one dollar, i.e. $h=1$. For contradiction, let us also assume that the overshipment cost is less than the holding cost, e.g. $c_o = 0.7$. In this case, it will always be optimal to ship all of the food out from the food bank, since it is less costly to waste it than to hold it at the food bank. In practice, food bank warehouses are significantly larger than the storage spaces of the charitable agencies in the counties, so holding food at the warehouse is cheaper than wasting it. Hence, we will assume $h \leq c_o$.

As mentioned previously, the optimal solution (\ref{eq:Xstar}) of the stochastic single-period food-distribution problem is analogous to that of the newsvendor problem. Consider the newsvendor problem with an upper limit on order quantity, i.e, purchase power \citep{birge1997introduction}. The analogy is revealed by representing the basic elements of the newsvendor problem using the same notation that we used to describe the food-distribution problem: (i) $S$ represents the total pounds of food available (purchase power), (ii) let $c^{NV}_U = h$ denote the cost of underage for the newsvendor problem; (iii) let $c^{NV}_O = c_o - h$ denote the cost of overage for the newsvendor problem where the difference represents the net cost of having shipped too much food to the counties, (iv) $\tilde{X}$ denotes the order quantity which is the decision variable, and (v) $R \Delta = D$ denotes the demand, which is a random variable with c.d.f. $F_D (x)=Pr\{D \leq x\}= Pr\{R \leq \frac{x}{\Delta}\}=F_R(x/\Delta)$  and p.d.f. $f_D(x)=\frac{1}{\Delta}f_R(\frac{x}{\Delta})$ for each $x \in X$. As $c^{NV}_U$ increases, the second term in (\ref{eq:Xstar}) also increases, resulting in a higher total shipment $\tilde{X^*}$; but eventually we reach the supply-constrained solution in which $\tilde{X^*} = S \leq \Delta F_R^{-1} (\frac{c_U}{c_U+c_O})$.  This is expected because if it is too costly to hold unshipped food at the food bank, the optimal solution should prescribe shipping the entire food supply to the counties. Moreover, if the cost $c^{NV}_O$ of overshipment increases, the second term in (\ref{eq:Xstar}) decreases; and we eventually reach a capacity-constrained solution in which $\tilde{X^*} = \Delta F_R^{-1} (\frac{c_U}{c_U+c_O}) \leq S$.  This is also intuitive because if the cost of overshipment is too high, then the optimal solution should prescribe shipping less food from the food bank. The particular values for these cost coefficients in practice depend on the preferences of the decision maker and can be estimated using statistical learning methods \citep{james2013introduction} based on previous data or by trial-and-error. In an extreme case, if the decision maker would like to avoid any overage, the cost of overage could be set as $c_o = \infty$, resulting in $\tilde{X}^* = \min\{S,r_L\Delta\}$. Conversely, if the decision maker would like to avoid any holding cost, since $h \leq c_o$ by definition, we could set $h= c_o$, which would result in $\tilde{X}^* = \min\{S,r_U\Delta\}$. In our computational experiments, we will typically set $h=1$ and vary $c_o$ to explore the change in distribution policies. In the next section, we extend this approach to multiple periods, yielding a more realistic representation of food bank operations.

\emph{Why do we only penalize overshipments and not undershipments?} We have seen from Proposition \ref{prop:SPMX} and Remark \ref{obs:1} that if there exists at least one county where overshipment occurred, all other counties must have also received overshipments (there cannot be any undershipments). 
We can similarly show that if at least one county received an undershipment, all other counties also received undershipments. That is, if $X_j < \hat{R} D_j$ for at least one county $j$, then $X_i < \hat{R} D_i$ for all other counties as well. However, since by definition, we have $R \leq \frac{\hat{C_j}}{D_j}$, this implies that $X_j \leq \hat{C_j}$ for all $j$. Thus, in this case, the food has been distributed equitably per \eqref{con:equitySPM} and the realized capacities did not result in any waste. The only issue with this case is the missed opportunity to ship more food; however, this is already accounted for through the holding cost for the unshipped food in our objective function.

\subsection{Multi-Period Problem}\label{sec:multiperiod}

The multi-period stochastic food distribution problem considers a finite planning horizon of $T$ periods. The decisions in each period follow those described in Section \ref{sec:model}, except that any inventory remaining at the food bank at the end of the period can now be carried over as supply for later periods. In addition to stochastic capacities, we now also consider uncertainty in future supply (food donations). Decisions for the period $t$ are made after donations $S_t$ for period $t$ have been received, but future donations, $S_\tau$ for $\tau>t$ are not known with certainty. This reflects the real-world operations of food banks, which typically make shipment decisions after receiving food donations for a given period. However, they may not have exact information about the amounts of donations to be received in future periods.
The following random variables are now time-dependent:

\quad \quad $C_{jt}:$ Capacity of county $j$ at time $t$,

\quad \quad $\frac{C_{jt}}{D_j}:$ CD ratio of county $j$ at time $t$,

\quad \quad $R_{t} = \min_{j \in \mathcal{J}} {\frac{C_{jt}}{D_j}}:$ Minimum CD ratio over all counties in the service region at time $t$,

\quad \quad $S_{t}$: donations received at the food bank at the start of period $t$.

We assume that the given (deterministic) demand values for counties $D_j$ are not time-dependent, since these values are generally estimated using the poverty populations in each county \citep{USCensus} and are fairly stable over short periods of time, such as the planning horizons considered in this paper. 
We define the following decision variables:

\quad \quad $X_{jt}:$ Amount of shipment (pounds) to county $j$ at time $t$, where $\tilde{X_t} = \sum_{j=1}^n {X_{jt}}$,

\quad \quad $Y_{jt}:$ Overshipped amount (pounds) to county $j$ at time $t$, where $\tilde{Y_t} = \sum_{j=1}^n {Y_{jt}}$.

The value $\hat{C_{jt}}$ denotes the observed capacity of county $j$ at time $t$, and $\hat{R_t}$ the minimum observed CD ratio in period $t$. The sequence of events in period $t$ is similar to that in the single-period problem, but involves time-dependent quantities as follows:
\begin{enumerate}
    \item The starting inventory level $\alpha_t$ given by
     \begin{equation} \label{eq:AlphaDef}
    \alpha_t = \alpha_{t-1} + S_{t-1} - \tilde{X}_{t-1} \quad \quad t=2,\dots,T
    \end{equation}
   is observed, where the starting inventory level $\alpha_1$ at the start of period 1 is a known parameter.
    \item The supply replenishment of $S_t$ pounds is received, raising the available inventory to $\beta_t=\alpha_t+S_t$.
    \item A total shipment of $\tilde{X_t} = \sum_{j=1}^n {X_{jt}}$  is sent to counties such that each $X_{jt}=\tilde{X_t} \frac{D_j}{\Delta}$  for $j \in \mathcal{J}$ and $\tilde{X_t} \leq $ $\beta_t$.
    \item The capacities of the counties, $\hat{C_{jt}}$ and hence, the minimum CD ratio $\hat{R_{t}} = \min_{j \in \mathcal{J}} {\frac{\hat{C_{jt}}}{D_j}}$ are observed.
    \item The amount of overshipped food at time $t$, $\tilde{Y_t} = \sum_{j=1}^n {Y_{jt}}$, is determined considering the observed capacities of the counties, while preserving perfect equity such that $ \frac{X_{jt}-Y_{jt}}{\tilde{X_t}-\tilde{Y_t}} - \frac{D_{j}}{\Delta} = 0$ and $X_{jt} - {\hat{C_{jt}}} \leq Y_{jt}$  for $j \in \mathcal{J}$.
    \item The overshipment and holding costs are computed for the given values of $\tilde{X_t}$ and $\hat{C_{jt}}$. 
\end{enumerate}

We assume that the minimum CD ratios $\{R_t:t=1,\dots,T\}$ are independent random variables with known c.d.f. $F_R(r)$ for $r \in [ r_L,r_U]$, where $0 \leq r_L < r_U \leq \infty$ (we take $r \in [r_L,\infty)$ if $r_U=\infty)$, with expected value $E[R]<\infty$; and when the minimum CD ratio $R_t$ is a continuous random variable, it has a given p.d.f. $f_R(r)$ that satisfies hypotheses (i) and (ii) of Proposition 1.  The per-pound cost coefficients $c_o$ and $h$ are defined as in Section \ref{sec:singleperiod} and are not time-dependent. We also assume that the supply (donations) received in future periods $\tau > t$ are also independent random variables with given c.d.f. $F_{S_\tau}(s_\tau)$ and p.d.f. $f_{S_\tau}(s_\tau)$.

Let $\beta_t$ represent the available inventory at the food bank after the donations at period $t$ have been received, i.e. $\beta_t = \alpha_t + S_t$. At the start of period $t \geq 1$, we must decide how much food will be shipped to counties to achieve $Q_t(\beta_t)$, the minimum expected total cost incurred over the remaining decision periods $t, t+1, \dots, T$, given that the food bank has $\beta_t$ pounds of food available to ship at the start of the period. 
\begin{equation}
\begin{aligned}\label{eq:Qt1}
    Q_t(\beta_t)= \min_{0 \leq x \leq \beta_t} \{h(\beta_t - x) + g_t(x)  + \mathbb{E}_{S_{t+1}}[Q_{t+1}(\beta_{t+1}=\beta_t-x+S_{t+1})]\}
\end{aligned}
\end{equation}
where
\begin{equation} \label{eq:gt1}
    g_t(x)= c_o x -  c_o \int_{r_L}^{\frac{x}{\Delta}} r \Delta f_R(r) dr - c_o \int_{\frac{x}{\Delta}}^{r_U} x f_R(r) dr
\end{equation}
which follows from equation \eqref{con:Z-mansc} for $x \in [0,\beta_t]$ and $t=1,\dots,T$. The following discussions will still hold if $R_t$ is a discrete random variable with a probability mass function $f_R(r)$. 

We now consider what happens at the end of the planning horizon. Considering (\ref{eq:Qt1}) for the last period $t=T$, we set a maximum ending inventory level, $\Theta$, and penalize overages exceeding this inventory level by setting
\begin{equation} \label{eq:Ending1}
    Q_{T+1}(\beta_{T+1}) = h_{T+1} [\beta_{T+1}- \Theta]^+.
\end{equation}

The maximum ending inventory level, $\Theta$, represents the maximum inventory that the food bank would ideally want to have at the end of the planning horizon. This is based on the common practice of food banks wanting enough storage space to accommodate incoming donations, as they typically do not refuse donations. By definition, if $S_{t+1} > \Theta$, then the term $\beta_{T+1}- \Theta$ will always be positive regardless of the decision maker's past shipment decisions. We now make some observations about this model.

\begin{lemma} \label{obs:history}
For the multi-period stochastic food distribution problem, at a given period $t$, where $t \geq 2$, the decisions in periods $t,\dots,T$ depend on the past decisions \{$\tilde{X_\tau}: \tau=1,\dots,t-1$\} only through
\begin{equation} \label{eq:alpha1}
    \alpha_t = \alpha_1 + \sum_{\tau=1}^{t-1} S_\tau - \sum_{\tau=1}^{t-1} \tilde{X_\tau}.
\end{equation}
Therefore knowledge of the entire shipment history \{$\tilde{X_\tau}: \tau=1,\dots,t-1$\} is not needed.
\end{lemma}

Lemma \ref{obs:history} follows directly from (\ref{eq:AlphaDef}) and implies that if at the start of any period $t$ we know how much food we are going to ship, the starting inventory level for the next period follows a stochastic process depending only on the supply received at the beginning of the next period. 
This suggests that the multi-period food distribution problem can be optimally solved using dynamic programming (DP) recursing backwards from the end of the planning horizon, $T+1$. However, in order to formulate this problem as a DP, we need to discretize the sample and action space. For this purpose, we make the following assumptions:

\begin{enumerate}[wide, labelwidth=!, labelindent=0pt]
    \item Minimum CD ratios $R_t$ are represented as discrete random variables with a known probability mass function $f_R(r)$ defined in the sample space $\Psi$.
    \item Donations $S_t$ are represented as discrete random variables with known probability mass functions $f_S(s)$ defined in the sample space $\rho$.
    \item The state space $\beta_t$ and action space $\tilde{X}_t$ are discrete.
\end{enumerate}

We thus formulate the DP to optimally solve this problem as follows \citep{snyder2019fundamentals}:

\noindent \underline{\textbf{The Dynamic Program for Solving the Multi-Period Food Distribution Problem}}
\begin{enumerate}[wide, labelwidth=!, labelindent=0pt]
    \item Require $\alpha_1, \{D_j: j=1,\dots,n\}, h, c_o, h_{T+1}, \Theta$. 
    \item Calculate the terminal cost function using (\ref{eq:Ending1}) for all $\beta_{T+1}$ and set $t \leftarrow T$.
    \item Compute $Q_t(\beta_t)$ for all $\beta_t$ as follows:
    \begin{equation}
    \begin{aligned}   
         Q_t(\beta_t)= \min_{0 \leq x \leq \beta_t} \{h(\beta_t - x) +  c_o x -  c_o \Delta \sum_{r_L \leq r_t^i \leq \frac{x}{\Delta}} {r_t^i p_t^i} - c_o  \sum_{\frac{x}{\Delta} \leq r_t^i \leq r_U} x p_t^i \\ + \sum_{S_{t+1} \in \rho}{f_{S_{t+1}}(s) Q_{t+1}(\beta_t-x+S_{t+1})}\}   \end{aligned}
    \end{equation}
    and set $\tilde{X_t} = \text{argmin} \{ Q_t(\beta_t) \}$.
    \item Set $t \leftarrow t-1$. If $t=0$, STOP. Otherwise, go to 3.
    \item Output $\{\tilde{X_t}^*=\tilde{X_t}: t=1,\dots,T\}$ and $Z^* = Q_1(\alpha_1 + S_1)$.
\end{enumerate}


The minimum expected total cost for the entire time horizon is given by $Q_1(\alpha_1 + S_1)$, where $\alpha_1$ is the known starting inventory at time $t=1$. 
The DP requires calculating the cost functions and specifying actions to be taken for each period based on all possible starting inventory levels. Populating all the possible inventory levels in each period, even if we discretize levels frugally, requires considerable computational effort bringing with it the well-known curse of dimensionality. Food banks usually lack the resources to implement complex algorithms and seek solutions that are easy to understand and implement. 
For these reasons, in the next sections, we develop heuristics to provide solutions for the multi-period food distribution problem. 

\section{Heuristics}\label{sec:methods}

\subsection{A Myopic Heuristic} \label{sec:myopic}
Our first heuristic uses the structure of the optimal solution to the single-period problem, presented in Section \ref{sec:singleperiod}, to propose a Myopic Heuristic for the multi-period problem assuming that the c.d.f. $F_R(r)$ of the minimum CD ratios is known. The total amount of food shipped in period $t$ by the Myopic Heuristic is denoted by $\tilde{X}_t^{MH}$ and the associated shipment to county $j$ in period $t$ by $X_{jt}^{MH}$.  The shipment $\tilde{X}_t^{MH}$ from the food bank in period $t$ must be made using the available supply $\beta_t^{MH} = \alpha_t^{MH}+S_t$ of food available at the start of period $t$. The remaining parameters have the same definitions as in Section \ref{sec:multiperiod}. This algorithm, being myopic in nature, does not take future supply uncertainty into consideration, and can be stated as follows:

\noindent \underline{\textbf{Myopic Heuristic (MH)}}
\begin{enumerate}[label*=\arabic*.]
    \item Require $\alpha_1, \{D_j: j=1,\dots,n\}, h, c_o, h_{T+1}, \Theta$.
    \item For each period $t=1,\dots,T$:
    \begin{enumerate}[label*=\arabic*.]
        \item \textbf{\textit{Observe current available inventory.}} Calculate the current available inventory as the following
        \begin{equation}
            \beta_t^{MH} = \left\{ 
                \begin{array}{ c l }
                \alpha_1 + S_1 & \quad \textrm{if } t=1 \\
                \beta_{t-1}^{MH}-\tilde{X}_{t-1}^{MH}+S_{t} & \quad \textrm{if } t > 1
                \end{array}
            \right.
        \end{equation}
        \item \textbf{\textit{Calculate shipments.}} Calculate the total shipment from the food bank \\ $\tilde{X}_t^{MH}=\min\{\beta_t^{MH},\Delta F_R^{-1}(\frac{h}{c_o})\}$ and set $X_{jt}^{MH}=\frac{\tilde{X}_t^{MH} D_j}{\Delta}$.
    \end{enumerate}
    \item \textbf{\textit{Objective Function Calculation.}} Calculate the objective function value as 
    \vspace{-10pt}
    \begin{equation} \label{eq:MHZ*}
    Z^{MH} = \sum_{t = 1}^{T} (h(\beta_t - \tilde{X}_t^{MH})+g_t(\tilde{X}_t^{MH})) + Q_{T+1}(\beta_{T+1}).
    \end{equation}
\end{enumerate}

The MH approaches the multi-period problem by solving a sequence of single-period problems, independent of future periods. Its obvious drawback is lack of any visibility of the changes in supply as well as costs that may be incurred if there is insufficient food to ship at a future period or too much supply remaining for a future period. We now develop structural properties and derive bounds on the deviation of the Myopic Heuristic from the optimal solution.

\subsection{Structural Properties and Bounds} \label{sec:structural}

Recall that $\beta_{t}$ represents the amount of food available for distribution at the start of period $t$ after donations for that period $S_t$ have been received. Lemma \ref{obs:history} implies the following corollary: 



\begin{corollary} \label{cor:indep}
For the multi-period food distribution problem, given $\beta_{t+1}$, decisions made in future periods $\tilde{X}_\tau: \tau=t+1,\dots,T$ are independent of the decisions $\tilde{X}_t$ made in period $t$.
\end{corollary}




Let $U_t(\beta_t,\tilde{X}_t)$ represent the objective function value of the single-period model presented in Section \ref{sec:singleperiod} that is optimally solved by the Myopic Heuristic in period $t$ and $V_t(\beta_t,\tilde{X}_t)$ the objective function value  of the multi-period problem over periods $t,\dots,T$. Let $\beta_t^{MH}$ denote the inventory available at the food bank for shipment at the start of period $t$ if the Myopic policy is implemented in periods $1,\dots,t-1$ and $\beta_t^*$ the inventory available for shipping at the start of period $t$ in an optimal solution to the multi-period problem. Then,
\begin{equation} \label{eq:Ut}
    U_t(\beta_t^{MH},\tilde{X}_t) = h(\beta_t^{MH}-\tilde{X}_t) + \mathbbm{E_\textbf{C}}[c_o\tilde{Y}_t]\quad \forall t,
\end{equation}
\begin{equation}
\begin{aligned} \label{eq:Vt}
V_t(\beta_t^{*},\tilde{X}_t) = 
     h(\beta_t^{*} - \tilde{X}_t) + \mathbbm{E_\textbf{C}}[c_o\tilde{Y}_t] + \mathbbm{E_{S_{t+1}}[V_{t+1}(\beta_{t+1}^*,\tilde{X}_{t+1})]}
\end{aligned}
\end{equation}
for $t=1,\dots,T$ where $\beta_{t+1}=\beta_t-\tilde{X}_t+S_{t+1}$ and 
\begin{equation} \label{eq:VEnd}
    V_{T+1}(\beta_{T+1}^{*},\tilde{X}_{T+1}) =V_{T+1}(\beta_{T+1}^{*})= h_{T+1} [\beta_{T+1}- \Theta]^+
\end{equation}
since no shipments are made in period $T+1$. The objective function value for the Myopic Heuristic in period $t$, which corresponds to the optimal solution to the single period problem with the same inputs, is:
 \begin{equation} \label{eq:Ht}
     H_t(\beta_t^{MH}) = \min_{\tilde{X}_t \leq \beta_t^{MH}} \{U_t(\beta_t^{MH},\tilde{X}_t)\}
 \end{equation}
On the other hand, the optimal solution in any period $t$ will minimize the following cost-to-go function for the periods $t,…,T$:
\begin{equation}\label{eq:Ht}
\begin{aligned}
B_t(\beta_t^*)
&= \min_{\tilde{X}_t \leq \beta_t^{*}}
\left\{V_t(\beta_t^{*},\tilde{X}_t)\right\} \\
&= \min_{\tilde{X}_t \leq \beta_t^{*}}
\left\{
U_t(\beta_t^{*},\tilde{X}_t)
+ \mathbb{E}_{S_{t+1}}
\left[
V_{t+1}(\beta_{t+1}^*,\tilde{X}_{t+1})
\right]
\right\} \\
&= \min_{\tilde{X}_t \leq \beta_t^{*}}
\left\{
U_t(\beta_t^{*},\tilde{X}_t)
+ \mathbb{E}_{S_{t+1}}
\left[
U_{t+1}(\beta_{t+1}^{*},\tilde{X}_{t+1}^{*})
+ \mathbb{E}_{S_{t+2}}
\left[
U_{t+2}(\beta_{t+2}^{*},\tilde{X}_{t+2}^{*})
+ \dots
+ \mathbb{E}_{S_T}
\left[
U_T(\beta_T^{*},\tilde{X}_T^{*})
+ \mathbb{E}_{S_{T+1}}
\left[
V_{T+1}(\beta_{T+1}^{*})
\right]
\right]
\right]
\right]
\right\}.
\end{aligned}
\end{equation}


Let 
\begin{equation} \label{eq:XtMH}
    \tilde{X}_t^{MH} = \argmin_{\tilde{X}_t \leq \beta_t^{MH}} \{U_t(\beta_t^{MH},\tilde{X}_t)\}
\end{equation}
\begin{equation} \label{eq:Xt*}
    \tilde{X}_t^* = \argmin_{\tilde{X}_t \leq \beta_t^{*}} \{V_t(\beta_t^{*},\tilde{X}_t)\}
\end{equation}
denote the minimizers of the functions $U_t(\beta_t^{MH},\tilde{X}_t^{MH})$ and $V_t(\beta_t^*,\tilde{X}_t^{*})$ subject to all constraints in the respective models. In other words, $\tilde{X}_t^{MH}$ minimizes the Myopic objective in period $t$ and $\tilde{X}_t^*$ the cost-to-go function for periods $t,\dots,T$. Note that $V_t(\beta_t^*,\tilde{X}_t^{*})$ represents the objective function of the multi-period problem in period $t$ and $B_t(\beta_t^*)$ the optimal objective function value of the multi-period problem in period $t$. Similarly, $U_t(\beta_t^{MH},\tilde{X}_t^{MH})$ represents the objective function of the Myopic Heuristic in period $t$ and $H_t(\beta_t^{MH})$ the optimal value of the Myopic objective function in period $t$. Proposition \ref{prop:MHLB} states that for any given value of $\beta_t$ the arg min of $U_t(\tilde{X}_t,\beta_t)$ never exceeds that of $V_t(\tilde{X}_t,\beta_t)$. The proof is given in Appendix \ref{Appendix:3}.

\begin{proposition}\label{prop:MHLB}
Define $\tilde{X}_t^{MH}$ and $\tilde{X}_t^*$ as in (\ref{eq:XtMH}) and (\ref{eq:Xt*}) and let $\beta_t^{MH}=\beta_t^*=\beta_t$ for period $t$ for $t=1,\dots,T$. Then,
the total distribution by the Myopic Heuristic in period $t$ is a lower bound on the optimal total distribution in that period for the multi-period problem, i.e., 
    \begin{equation}\label{eq:Prop1i}
        \tilde{X}_t^{MH} \leq \tilde{X}_t^*.
    \end{equation}
Further, if  all holding costs are zero, i.e. $h=h_{T+1}=0$, we will have $\tilde{X}_t^{MH} = \tilde{X}_t^*$.
\end{proposition}

Proposition \ref{prop:MHLB} implies that the MH will never ship more than the optimal solution as long as holding costs are positive. This result is intuitive because the myopic solution does not consider future inventory holding costs at the food bank and hence gives more conservative solutions (in terms of shipment levels) than the optimal solution. However, if $h=h_{T+1}=0$, we will have $\tilde{X}_t=\tilde{X}_t^{MH}$ and the myopic solutions will be optimal for the multi-period problem. This is reasonable as in an environment where there is no holding cost, the only remaining cost is that associated with overshipping. In this case, the food bank would choose to ship any value between $[0,\min\{r_L\Delta,\beta_t\}]$, which would result in a zero overshipment cost. This solution would be optimal for both the myopic and the multi-period problems since there is no holding cost. Thus, we would expect the MH to perform worse in high holding cost scenarios. We explore this further in Section \ref{sec:experiments}. Finally, since we have a minimization problem, the objective function value $V_t(\beta_t,\tilde{X}_t^{MH})$ generated by the Myopic Heuristic provides an upper bound on the optimal objective function value $V_t(\beta_t,\tilde{X}_t^{*})$.

We now develop an upper bound on the optimal shipment in period $t$, yielding a lower bound on the optimal objective function value.

\begin{proposition}\label{prop:MLMUB}
Consider the cost-to-go function defined in (\ref{eq:Vt}). 
If we set $c_o=0$ such that 
\begin{equation} \label{eq:Vt_UL}
    V_t^L(\beta_t,\tilde{X}_t) = 
     h(\beta_t - \tilde{X}_t) + \mathbbm{E_{S_{t+1}}}[V_{t+1}(\beta_{t+1},\tilde{X}_t)]
\end{equation}
we have the following results:
\begin{enumerate}
    \item The minimizer of $V_t^L(\beta_t,\tilde{X}_t)$ is given by
    \begin{equation}\label{eq:XU}
       \argmin_{0 \leq \tilde{X}_t \leq \beta_t}  V_t^L(\beta_t,\tilde{X}_t) = \tilde{X}_t^U = \beta_t
    \end{equation}
    \item We have $\tilde{X}_t^* \leq \tilde{X}_t^U$, i.e. the minimizer of $V_t^L(\beta_t,\tilde{X}_t)$ provides an upper bound on the optimal total shipment in period $t$.
    \item Finally, the cost-to-go function \eqref{eq:Vt_UL} provides a lower bound on the optimal cost-to-go function, i.e., $V_t^L(\beta_t,\tilde{X}_t^U) \leq V_t^*(\beta_t,\tilde{X}_t^*)$.
\end{enumerate}
\end{proposition}

With no overage cost, the solution that ships the entire available food supply each period, i.e. $\tilde{X}_t^U = \beta_t$ results in a zero-cost solution, proving part 1. Since we also have $\tilde{X}_t \leq \beta_t$ as the supply constraint at each period, part 2 is also true. Finally, it is straightforward to see that a zero-cost solution will provide a lower bound on the optimal solution $V_t^*(\beta_t,\tilde{X}_t^*)$.

Proposition \ref{prop:MHOpt} defines an additional situation where the MH will be optimal.

\begin{proposition}\label{prop:MHOpt}
Assume that optimal decisions have been taken prior to time $t$ and the current available inventory is $\beta_t^*$. The Myopic Heuristic provides the optimal solution, i.e. $\tilde{X}_t^{MH}=\tilde{X}_t^{*}$ if the following condition holds:
\begin{equation}
    r_L \geq \frac{\beta_t^*}{\Delta}.
\end{equation}
\end{proposition}

\begin{proof}
We have already established that $\tilde{X}_t^{MH} \leq \tilde{X}_t^* \leq \tilde{X}_t^U$ (Propositions \ref{prop:MHLB} and \ref{prop:MLMUB}) where $\tilde{X}_t^{MH} = \min\left\{\beta_t^*,\Delta F_R^{-1}\left(\frac{h}{c_o}\right)\right\}$,
and $\tilde{X}_t^{U} = \beta_t^*$. It is straightforward to see that if the lower and upper bounds are equal, i.e., $\tilde{X}_t^{MH} = \tilde{X}_t^U$, then the solution generated by the MH will be optimal, i.e., $\tilde{X}_t^{MH} = \tilde{X}_t^U = \tilde{X}_t^*$. This would be achieved if: $\Delta F_R^{-1}\left(\frac{h}{c_o}\right) \geq \beta_t^*$, or, $F_R^{-1}\left(\frac{h}{c_o}\right) \geq \frac{\beta_t^*}{\Delta}$.
Now, if the condition given in Proposition \ref{prop:MHOpt} holds, i.e., $r_L \geq \frac{\beta_t^*}{\Delta}$, this implies that:
\begin{equation}
    F_R^{-1}\left(\frac{h}{c_o}\right) \geq r_L \geq \frac{\beta_t^*}{\Delta},
\end{equation}
The proof is completed.
\end{proof}

In practical terms, this proposition states that the MH solution will be optimal if the per-person fair share of our available inventory $\frac{\beta_t^*}{\Delta}$ does not exceed the lowest possible minimum CD ratio, implying that the problem is not constrained by the capacities but by available supply. In such a scenario, distributing all available inventory cannot lead to overshipment and gives zero holding cost, representing the best possible solution. 

To summarize our findings, let $\beta_t^*$ denote the on-hand inventory in any time period $t$ assuming that optimal shipment decisions have been made in periods $1,\dots,t-1$. We then have 
\begin{equation} \label{eq:boundsonX}
   \tilde{X}_t^{MH} \leq \tilde{X}_t^* \leq \tilde{X}_t^U,
\end{equation}
where
\begin{equation} \label{eq:sum2}
    \tilde{X}_t^{MH} = \min\left\{\beta_t^*,\Delta F_R^{-1}\left(\frac{h}{c_o}\right)\right\} \quad \text{and} \quad \tilde{X}_t^{U} = \beta_t^*.
\end{equation}
We have also shown that
\begin{equation}\label{eq:sum4}
    V_t^L(\beta_t,\tilde{X}_t^U) \leq V_t(\beta_t,\tilde{X}_t^*) \leq V_t(\beta_t,\tilde{X}_t^{MH}).
\end{equation}

Further, we have shown that the Myopic Heuristic provides optimal solutions to the multi-period problem if either of the following conditions hold:
\begin{enumerate}[label=\alph*.]
    \item  All holding costs are zero, i.e. $h=h_{T+1}= 0$, or
    \item  $r_L \geq \frac{\beta_t^*}{\Delta}$ for every $t$.
\end{enumerate}
However, neither of these conditions is realistic for food bank settings. Zero holding cost implies there is no cost to the food bank for holding the food indefinitely. However, holding food for long periods results in food spoilage and waste. Second, as we will later see in Section \ref{sec:experiments}, many food bank supply chains have portions of their service areas with insufficient capacity to serve the population in need. In Section \ref{sec:lookahead}, we develop a heuristic that considers the expected future supply to plan shipments for two periods at a time. Then, in Section \ref{sec:bound}, we introduce a heuristic that ships the average of $\tilde{X}_t^{MH}$ and $\tilde{X}_t^U$ in each period, providing a less conservative solution than the MH.

\subsection{The One-period Look-ahead Heuristic} \label{sec:lookahead}
In this section, we develop a heuristic which, at time $t$, considers the expected supply in the next period $t+1$, and adjusts shipments to ship more food if no supply shortfall is expected in the next period to reduce holding costs and food spoilage. We use the superscript ``LAH" to denote a solution obtained by the One-period Look-ahead Heuristic, referred to as the Look-Ahead solution. The differences between the presented heuristics are in the \textit{``Calculate shipments"} stage. Hence, to obtain the LAH, the ``MH" superscript in the Myopic Heuristic (Section \ref{sec:myopic}) is replaced by ``LAH" and Step 2.2 of the MH is replaced by the following:

\noindent \underline{\textit{\textbf{Shipment calculation step for the One-period Look-ahead Heuristic (LAH):}}}\\
\textbf{2.2.a} \textit{\textbf{Calculate myopic shipments.}} Calculate the total myopic shipment \\ $\tilde{X}_t^{LAH}=\min\{\beta_t^{LAH},\Delta F_R^{-1}(\frac{h}{c_o})\}$.\\
\textbf{2.2.b} \textit{\textbf{Check next period's supply.}} Calculate the expected surplus for the next period as:
\begin{equation}
\Phi_{t+1} =
\begin{cases}
\max\left\{
0,\,
\beta_t^{LAH} - \tilde{X}_t^{LAH}
+ \mathbb{E}_{S_{t+1}}[s_{t+1}]
- \Delta F_R^{-1}\left(\frac{h}{c_o}\right)
\right\},
& t = 1,\ldots,T-1, \\[6pt]
\max\left\{
0,\,
\beta_t^{LAH} - \tilde{X}_t^{LAH}
+ \mathbb{E}_{S_{t+1}}[s_{t+1}]
- \Theta
\right\},
& t = T.
\end{cases}
\end{equation}

If $\Phi_{t+1}>0$, increase the shipment in period $t$ as follows:
\begin{equation}
    \tilde{X}_t^{LAH}= \min\left\{\beta_t^{LAH}, \min\left\{\beta_t^{LAH},\Delta F_R^{-1}\left(\frac{h}{c_o}\right)\right\}+\Phi_{t+1}\right\}
\end{equation}
Finally, set $X_{jt}^{LAH}=\frac{\tilde{X}_t^{LAH} D_j}{\Delta}$.

The improvement of the LAH over MH is captured in Step 2.2.b. We have shown that the MH provides a lower bound on the shipments and may result in overages in future periods. In the LAH, we compute the expected surplus for the next period, assuming myopic shipments in the current period. Since supply is stochastic, this surplus calculation is based on an expected value of the supply variable. If making myopic shipments in the current period results in an expected surplus, and hence an overage cost in the next period, we increase shipments in the current period. Intuitively, in a situation where in a given period $t$, the food bank expects to receive a bulk donation in the next period $t+1$, the food bank may distribute more food in period $t$ in order not to be limited by distribution capacity in the next period. 

\subsection{The Bounding Heuristic} \label{sec:bound}
We introduce a heuristic that ships the average of the upper and lower bounds on shipments generated in Section \ref{sec:structural} for each period. This approach is motivated by the work of \citet{shang2003newsvendor} for multi-echelon serial production-inventory systems, and does not require any knowledge of future periods. The superscript ``BH" denotes a solution obtained by the Bounding Heuristic, referred to as the Bounding solution. In order to obtain the BH, the ``MH" superscript in the MH in Section \ref{sec:myopic} is replaced by ``BH" and Step 2.2 of the MH is replaced by the following:

\noindent \underline{\textit{\textbf{Shipment calculation step for the Bounding Heuristic (BH):}}}\\
\textbf{2.2.a} \textit{\textbf{Calculate bounds.}} Calculate the lower and upper bounds on the optimal shipment as $\tilde{X}_t^{L} = \min\{\beta_t^{BH},\Delta F_R^{-1}(\frac{h}{c_o})\}$, and $\tilde{X}_t^{U} = \beta_t^{BH}$ respectively.\\
\textbf{2.2.b} \textit{\textbf{Calculate shipments.}} Calculate the total shipment as
\begin{equation}\label{eq:BH}
    \tilde{X}_t^{BH} = \frac{\tilde{X}_t^{L}+\tilde{X}_t^{U}}{2}
\end{equation}
and set $X_{jt}^{BH}=\frac{\tilde{X}_t^{BH} D_j}{\Delta}$.

The Bounding Heuristic uses the average of the lower and upper bounds on the total optimal shipment and hence provides a less conservative policy than the Myopic heuristic, i.e. $\tilde{X}_t^{MH} \leq \tilde{X}_t^{BH}$. Hence we would expect BH to perform better in situations where MH results in overages in the following periods, similar to the LAH.

\section{Computational Experiments}\label{sec:experiments}
\subsection{Experimental Design}\label{sec:expdesign}
In this section, we first explore the computational time of the DP for the multi-period stochastic food distribution model for varying parameter settings. We then study the performance of the heuristics relative to the exact solution obtained from the DP. We design our experiments based on real-world data from our collaborator, FBCENC. 


We use fiscal year 2018-19 data on weekly distribution decisions and focus on dry goods (shelf-stable products such as pasta and canned goods) that constitute the largest proportion of food received and distributed by FBCENC. 
For the numerical study, we assume that the probability distributions of the minimum CD ratio $R_t$ and the weekly supply $S_t$ are both stationary over the period considered in this paper. Since FBCENC does not have data about the capacities of their counties, due to these being determined by the autonomous agencies in each county, we assume that the $90^{th}$ percentile $C_j^{90}$ of the weekly shipments to a county $j$ represents a reasonable estimate of that county's mean capacity \citep{orgut2016modeling}. We use the $90^{th}$ percentile rather than the maximum value since, based on our discussions with FBCENC, there are some weeks during the year, e.g. during the holiday season, where agencies receive much more food than their regular capacities can handle, causing them to temporarily acquire additional capacity in those extreme circumstances. Acquiring additional temporary capacity is expensive for charitable agencies operating on low budgets and is only done in rare circumstances. Hence, the maximum does not accurately reflect a county's capacity. We estimate demand using the average of the 2018 and 2019 poverty population estimates by the Census Bureau \citep{USCensus} for each county. Following Feeding America's guidelines \citep{feedingamericaMeals}, we assume that a person's need is approximately 3.6 pounds per day (25.2 pounds/week), so we multiply this number by the poverty population counts to estimate county $j$'s demand $D_j$. 
We assume that each county's CD ratio follows an Exponential distribution with rate $D_j/C_j^{90}$, similar to \citet{sengul2017modeling}. Since the Exponential distribution is closed under scaling, we obtain $R_t \sim \text{Exp}\left(\sum_{j \in J}{\frac{D_j}{C_j^{90}}}\right)$ = \text{Exp}(365.01).

We estimate the probability distribution of the weekly supply $S_t$ from FBCENC's weekly dry goods donation data during the 2018-19 fiscal year, a histogram of which is presented in Section~\ref{Appendix:HistDonation} of the supplementary material. Using the Seaborn\textsuperscript{TM} and Fitter\textsuperscript{TM} libraries in Python\textsuperscript{TM}, we find that a Gamma Distribution with shape parameter $1.9$, scale parameter $335276.48$, and location parameter $175150.42$ provides the best fit (lowest sum of squared error and BIC score). 


To implement the DP we discretize the probability distributions of $S_t$ and $R_t$ using the discrete approximation of \citet{miller1983discrete}, which is based on matching the first two moments of the probability distribution with the Gaussian Quadrature values tabulated in \citet{stroud1966gaussian}. 
For the minimum CD ratio $R_t$ we convert the tabulated values for the Exponential(1) distribution in \citet{stroud1966gaussian} to a $\text{Exponential}(365.01)$ distribution. We consider four discretization schemes by setting the number $Q_R$ of node-probability pairs for $R_t$ to values in the set $\{2,3,4,5\}$. For supply $S_t$ we use the Statmod\textsuperscript{TM} package in the R Programming Language\textsuperscript{TM} to generate the Gaussian Quadrature nodes and weights for the fitted Gamma distribution, and consider the same four discretization schemes as for $R_t$, i.e. $Q_S \in \{2,3,4,5\}$. The resulting nodes and weights for varying discretization schemes are given in Section~\ref{Appendix:5} of the supplementary material.

Solving the DP also requires discretization of the continuous state and action spaces \citep{cai2009dynamic}. We define a parameter $\lambda$ representing the granularity of discretization and explore its effect on the computation time. A lower $\lambda$ value should provide a better representation of the actual problem at the cost of higher computation time. We use the same granularity value $\lambda$ for both the state and action spaces.

To discretize the state space, consider the histogram of weekly donations to FBCENC (supplementary material \ref{Appendix:HistDonation}). Recall that the state space is defined by the possible values of the on-hand inventory $\beta_t$ at the start of period $t$ after donations for that period $S_t$ have been received. The amount of supply received in each period must be in the range $[\min_{\tau \in T}S_{\tau},\max_{\tau \in T}S_{\tau}]$. We round the lower bound down and the upper bound up to the nearest multiples of $\lambda$. Denoting these limits by $S^{min}$ and $S^{max}$, we can represent the state space for period $t$ as $[S^{min}, tS^{max}]$. Thus in any period $t$, the available inventory cannot be lower than the minimum donation amount $S^{min}$, and cannot exceed the case where the maximum donation amount of $S^{max}$ has been received in all periods up to the current period and no food has been shipped ($t*S^{max}$). 

Since we discretize the state and action spaces, we use piecewise linear interpolation \citep{cai2009dynamic} to estimate the value function when the realized state for a period $\beta_t$ is not captured in the state space. To do this, we use the value function $v^i$ associated with the highest possible available inventory value that is less than $\beta_t$ ($\beta_t^{i}$), and the value function $v^{i+1}$ associated with the next lowest possible available inventory ($\beta_t^{i+1}$). Hence, we approximate $V_t(\beta_t,\tilde{X}_t)$ by:
\begin{equation}
    V_t(\beta_t,\tilde{X}_t) = v_i + \frac{v^{i+1}-v^i}{\beta_t^{i+1}-\beta_t^{i}}(\beta_t-\beta_t^{i}).
\end{equation}
for  $\beta_t^{i} \leq \beta_t \leq \beta_t^{i+1}$.

Our experiments vary the length $T$ of the time horizon, overage cost $c_o$, terminal holding penalty $h_{T+1}$, the discretization scheme for $R_t$ ($Q_R$), the discretization scheme for $S_t$ ($Q_S$), the granularity of discretization for the state and action spaces $\lambda$, and ending target inventory $\Theta$.  Table \ref{tab:ExpDesign} summarizes our experimental design settings. We assume a starting inventory of $\alpha_1 = 0$.

\begin{table}[h]
    \small
    \centering
    \caption{Experimental Design Settings.}
    \begin{tabular}{ccc}\hline
    Parameter &	Settings\\ \hline
    $Q_R$ & \{2, 3, 4, 5\}\\
    $Q_S$ & \{2, 3, 4, 5\}\\
    $\lambda$ & \{1000,5000,10000,50000,100000\} pounds\\
    $T$&  \{1, 4, 8, 12, 16, 20\} weeks\\
    $c_o$ & \{1, 2, 3, 4, 5\} \\
    $h_{T+1}$ & \{1/3, 0.5, 1, 2, 3\}\\
    $h$ & 1\\
    $\Theta$ & \{$E[S]/3$, $E[S]/2$, $E[S]$, $2E[S]$, $3E[S]$\} pounds\\
    $R_t$ & $\omega$\emph{Exp}(365.01) where $\omega \in \{1/3, 1/2, 1, 2, 3\}$\\
    $S_t$ & \emph{Gamma}(1.9,335276.48,175150.42)\\\hline
    \end{tabular}
    \label{tab:ExpDesign}
\end{table}
    
\subsection{Numerical Results}\label{sec:numresults}
We first examine the relationship between problem size and the computation time and solution quality of the DP. We then compare the performance of the heuristics to the solution obtained via the DP under varying model parameters. 
All analyses were done on a 64-bit Intel Core i7 machine with 3.6 GHz processor and 64GB memory.

\subsubsection{DP Evaluation and Performance}\label{sec:DPComp}

In this section, we evaluate the computational performance of the DP for varying $\lambda$, $T$, $Q_R$, and $Q_S$. For the base case settings, we use $T=12$ weeks, $c_o=2$, $h=1$, $h_{T+1}=1$, $Q_R=4$, $Q_S=4$, $\lambda=5000$ pounds and $\Theta=E[S]$.

The size of the state space for a problem instance with $T$ periods can be calculated as $\sum_{t=1}^{T+1} \ceil{\frac{(t S^{max}-S^{min})}{\lambda}}$. Further, the number of actions to evaluate at each time period $t$ would be $\sum_{i=0}^{\ceil{t S^{max}/\lambda}}i$. Hence, the total number of actions to evaluate for the entire planning horizon is equal to 
\begin{align}\label{eq:ActSpace}
\begin{split}
    \sum_{t=1}^{T+1} \sum_{i=0}^{\ceil{t S^{max}/\lambda}}i \ceil{\frac{(t S^{max}-S^{min})}{\lambda}} = \sum_{t=1}^{T+1} \ceil{\frac{(t S^{max}-S^{min})}{\lambda}} \sum_{i=0}^{\ceil{t S^{max}/\lambda}}i \\= \sum_{t=1}^{T+1} \ceil{\frac{(t S^{max}-S^{min})}{\lambda}} \frac{\ceil{t S^{max}/\lambda} (\ceil{t S^{max}/\lambda}+1)}{2}
\end{split}
\end{align}
Hence, we would expect the computation time to increase polynomially as a function of $S^{max}, T$, and $\frac{1}{\lambda}$. 

\begin{figure}[h]
\begin{center}
\includegraphics[scale=0.25]{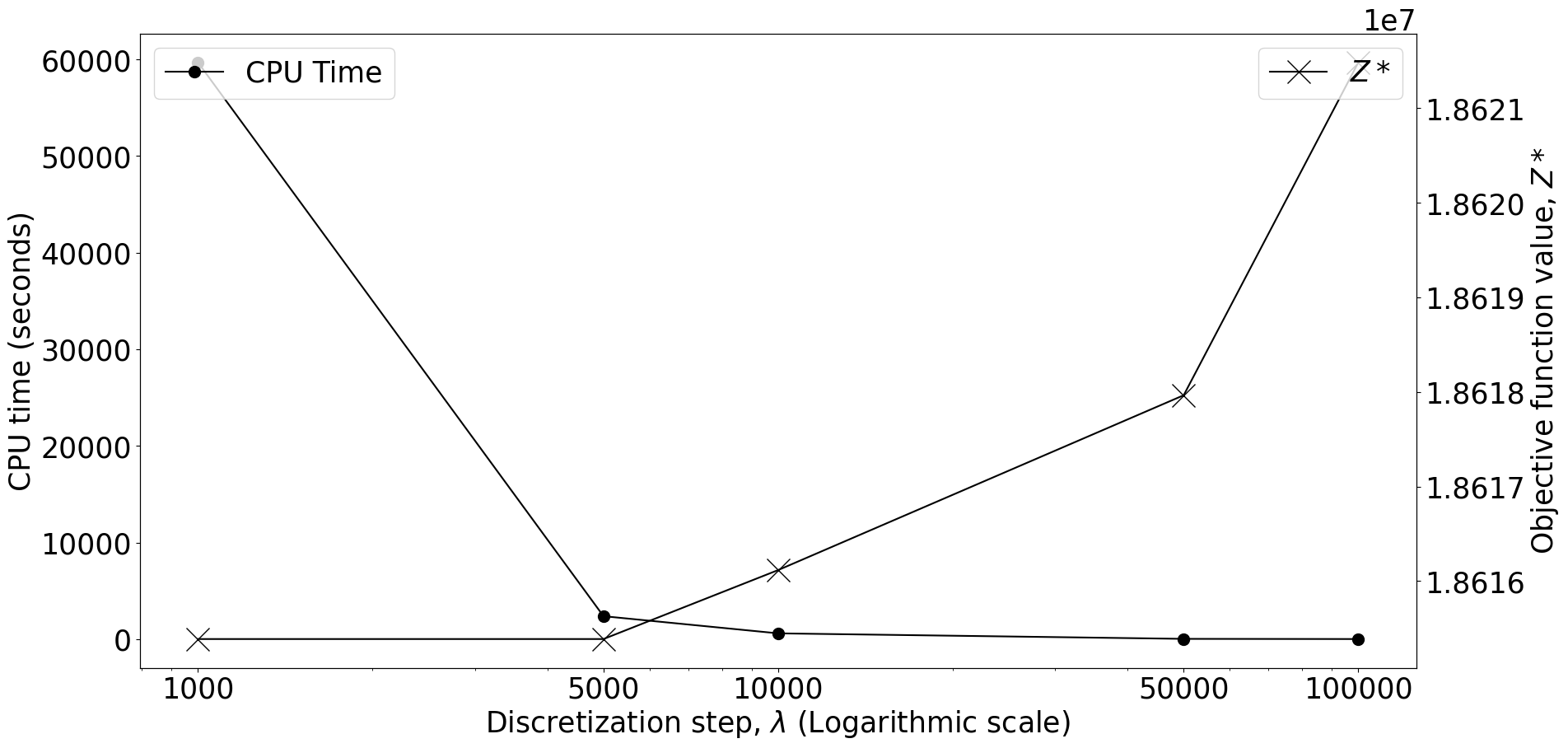}
\caption{The Impact of Changing the discretization granularity $\lambda$ on DP's Computational Performance and the Objective Function Value ($Z^*$).}
\label{fig:CompExp_CompTime_Lambda}
\end{center}
\vspace{-15pt}
\end{figure}

Figure \ref{fig:CompExp_CompTime_Lambda} illustrates the CPU time and the objective function value (expected cost) for the DP when we vary the discretization granularity $\lambda$. 
The computation time is very high when $\lambda=1000$ pounds but decreases rapidly as $\lambda$ increases, and does not change much after $\lambda=5000$ pounds. Since, based on (\ref{eq:ActSpace}), the number of actions to be evaluated is a function of $\frac{1}{\lambda^3}$, this is expected. The expected cost is non-decreasing in $\lambda$ with little change between $\lambda=1000$ and $\lambda=5000$ and increasing rapidly thereafter. A higher level of granularity in state and action spaces results in lower-cost decisions, as expected. We set $\lambda=5000$ in the following analyses as it achieves a good balance between solution quality and computation time. 

\begin{figure}[h]
\centering
\begin{tabular}{ll}
\subfloat[]{\includegraphics[width=3.2in]{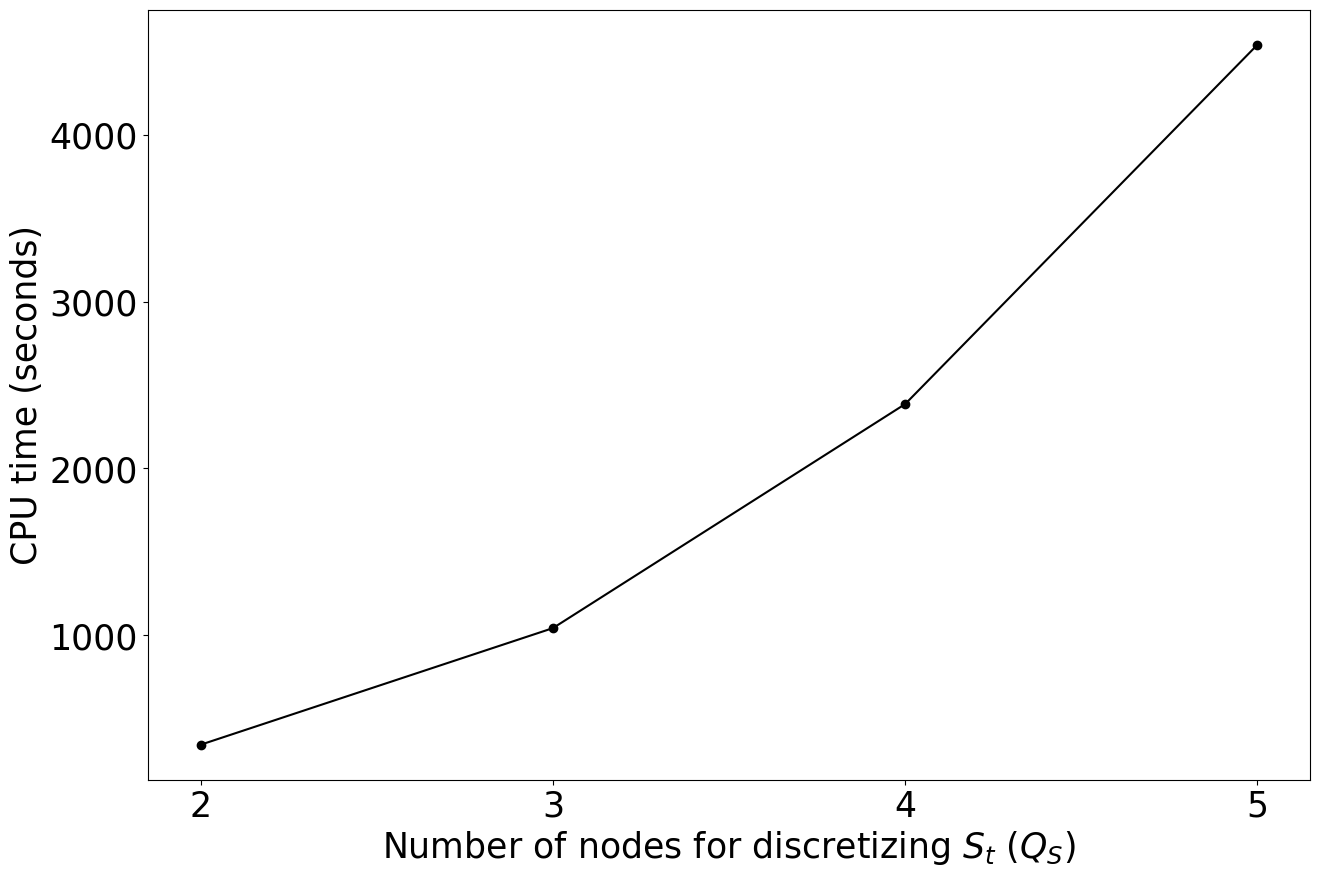}} &
\subfloat[]{\includegraphics[width=3.2in]{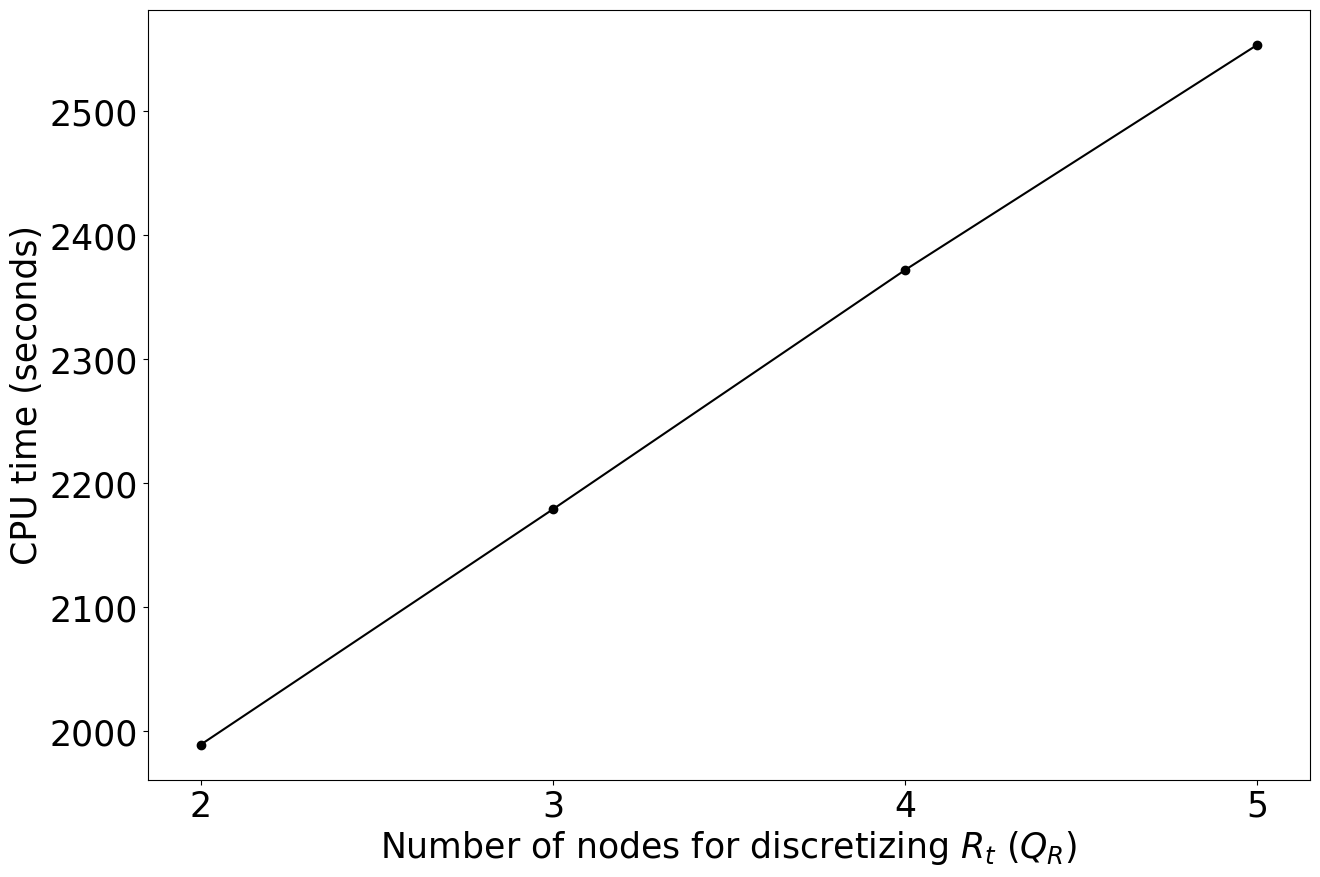}}
\end{tabular}
\caption{The Impact of Changing the discretization schemes $Q_S$ (a) and $Q_R$ (b) on the DP's Computational Performance.}
\label{fig:CompExp_CompTime_Q}
\end{figure}

As the number $Q_S$ of nodes for supply discretization increases, the range $[S^{min}, tS^{max}]$ becomes wider, increasing the computation time. Indeed, based on (\ref{eq:ActSpace}), we would expect the computation time to increase polynomially as $S^{max}$ increases and linearly as $S^{min}$ decreases. Note that the state and action spaces do not depend on the discretization of $R_t$, but a unit increase in  $Q_R$ requires one additional calculation for each state-action pair considered. Figure \ref{fig:CompExp_CompTime_Q} illustrates the CPU time of the DP when we vary $Q_S$ (Figure \ref{fig:CompExp_CompTime_Q}a) and $Q_R$ (Figure \ref{fig:CompExp_CompTime_Q}b).  
We see that computation time increases polynomially in $Q_S$ and almost linearly in $Q_R$. We set $Q_S=Q_R=4$ in the following analyses.  

\begin{figure}[h]
\begin{center}
\includegraphics[scale=0.25]{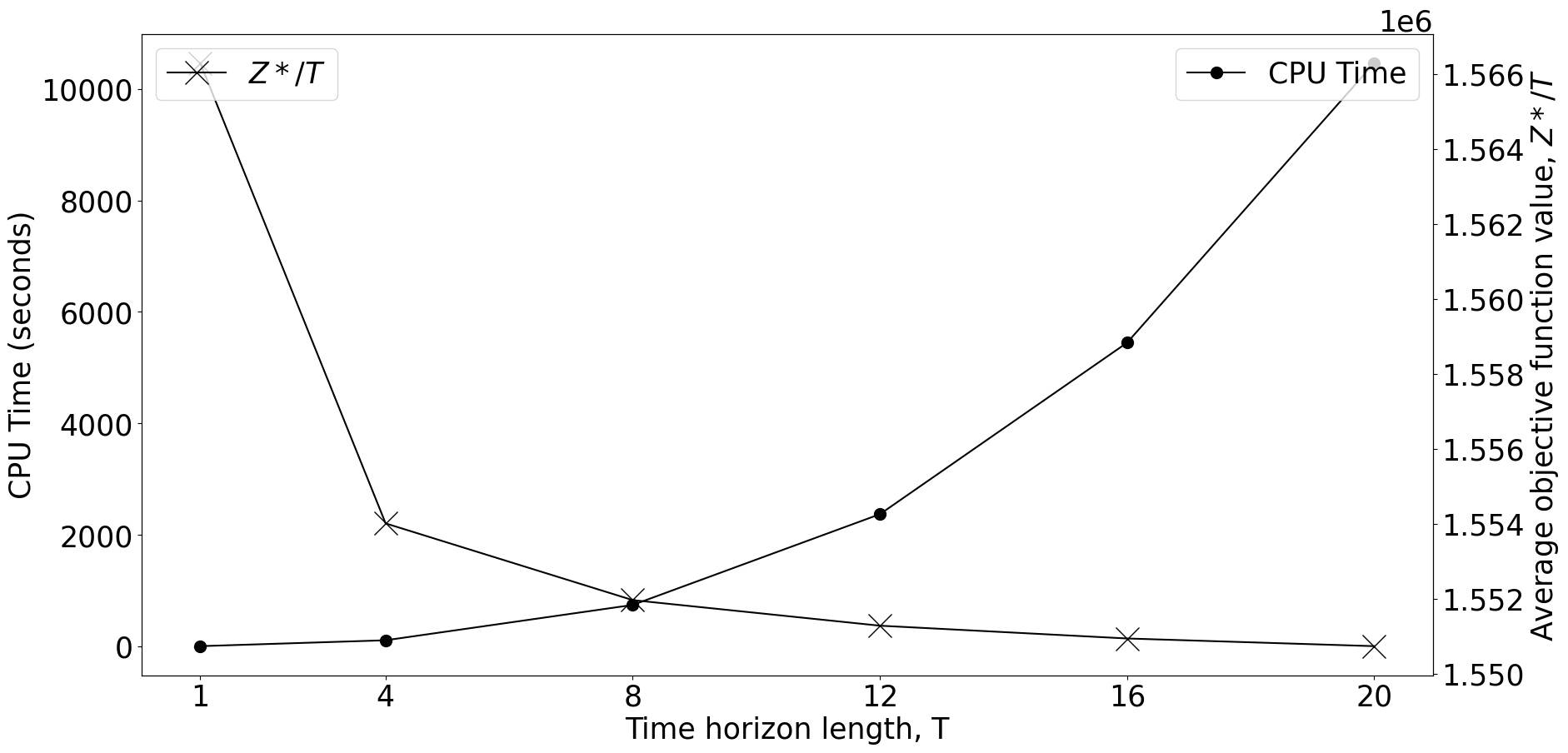}
\caption{The Impact of Changing the Time Horizon Length, $T$, on DP's Computational Performance and Objective Function Value.}
\label{fig:CompExp_CompTime_T}
\end{center}
\vspace{-15pt}
\end{figure}

Figure \ref{fig:CompExp_CompTime_T} illustrates the CPU time for DP and the average objective function value ($Z/T$) when we vary the time horizon length $T$. 
The average objective function value decreases as $T$ increases. This shows that planning for a longer horizon is beneficial for achieving a lower-cost solution and highlights the value of the multi-period model (over a single-period model). In fact, we see that the average $Z^*$ value decreases rapidly from $T=1$ to $T=8$ and does not change much after that. Hence, planning for a longer horizon is better, but the marginal benefit of increasing $T$ decreases after a certain point. As $T$ increases, the CPU times increase polynomially as expected. 
We set $T=12$ weeks in the following analyses as it results in a low-cost and low computation time solution, and aligns with the food banks' decision-making horizon, as it is common for food banks to make quarterly decisions. 

\subsubsection{Performance of the Heuristics }\label{sec:HeurComp}\mbox{}
After solving the DP using the parameter settings in Section \ref{sec:expdesign}, we generate 100 independent capacity and supply realizations to evaluate the performance of the heuristics. We illustrate our results using three metrics: ($i$) The objective function value (the total cost), ($ii$) The total pounds of food distributed, and ($iii$) The error ratio $\frac{Z^H-Z^{DP}}{Z^H} (100\%)$ where $Z^H$ denotes the objective function value for a given heuristic and $Z^{DP}$ the objective function value for the DP.

\begin{figure}[h]
\centering
\begin{tabular}{ll}
\subfloat[]{\includegraphics[width=3.2in]{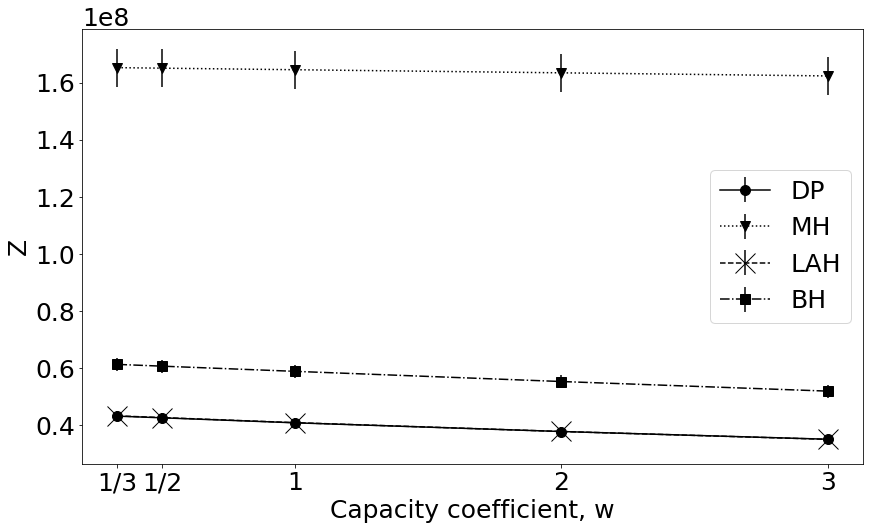}} &
\subfloat[]{\includegraphics[width=3.2in]{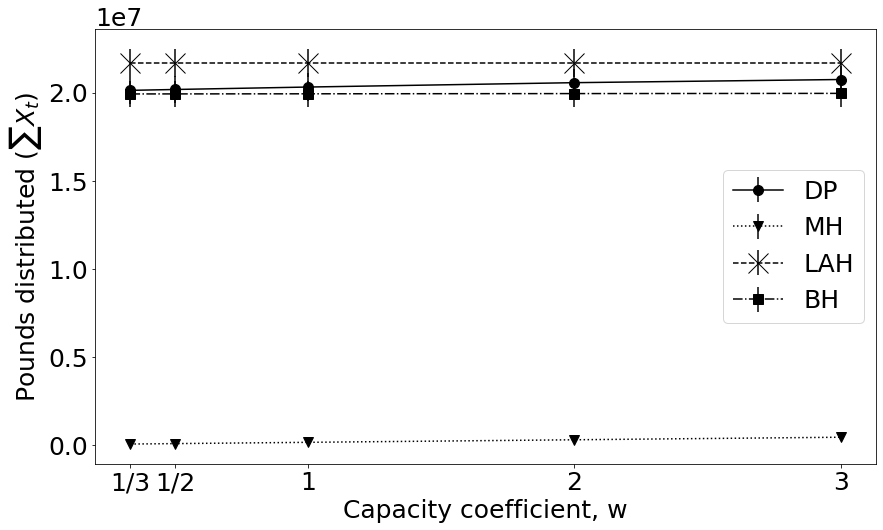}}
\end{tabular}
\caption{Results for (a) the objective function value and (b) the total food distribution over the time horizon for varying $w$.}
\label{fig:Heuristics_w}
\end{figure}

\noindent \emph{\textbf{Capacity Level, $w$.}} We vary the capacity in the service region by multiplying the probability distribution of the minimum CD ratio $R_t$ by a coefficient $w$. This allows us to test varying levels of capacity from a low-capacity scenario (when $w$ is small) to a high-capacity scenario (when $w$ is large). Figure \ref{fig:Heuristics_w} illustrates (a) the objective function value, and (b) the total food distribution over the time horizon for the DP, MH, LAH, and BH, while Table \ref{tab:ErrorRatio_w} presents the error ratios of the heuristics.

\begin{table}[htbp]
    \centering
    \small
    \setlength{\tabcolsep}{3pt}

    \caption{Statistics (mean (Avg), standard deviation (StD),
    minimum (Min), and maximum (Max)) for the error ratio
    for varying $w$.}
    \label{tab:ErrorRatio_w}

    \begin{adjustbox}{max width=\linewidth}
    \begin{tabular}{l|ccc|ccc|ccc|ccc|ccc}
    & \multicolumn{3}{c|}{$w=1/3$}
    & \multicolumn{3}{c|}{$w=1/2$}
    & \multicolumn{3}{c|}{$w=1$}
    & \multicolumn{3}{c|}{$w=2$}
    & \multicolumn{3}{c}{$w=3$} \\
    \hline
    & MH & LAH & BH
    & MH & LAH & BH
    & MH & LAH & BH
    & MH & LAH & BH
    & MH & LAH & BH \\
    \hline
    Avg &73.6\%&0.2\%&29.5\%
        &74\%&0.2\%&29.7\%
        &75\%&0.1\%&30.6\%
        &76.7\%&0\%&31.7\%
        &78.3\%&-0.2\%&32.5\%\\
    StD &2.4\%&0.4\%&1.6\%
        &2.4\%&0.4\%&1.6\%
        &2.4\%&0.5\%&1.7\%
        &2.5\%&0.7\%&2\%
        &2.7\%&0.9\%&2.7\%\\
    Min &65.8\%&-0.4\%&25.3\%
        &66\%&-0.6\%&25.3\%
        &66.4\%&-1.2\%&25.6\%
        &67.4\%&-2\%&26.2\%
        &68.4\%&-3.6\%&26.4\%\\
    Max &78.1\%&1.2\%&32.6\%
        &78.5\%&1.3\%&32.8\%
        &79.6\%&1.3\%&33.7\%
        &81.4\%&1.3\%&37.1\%
        &83.9\%&1.5\%&40.6\%\\
    \hline
    \end{tabular}
    \end{adjustbox}
\end{table}

Examining Figure \ref{fig:Heuristics_w}a and Table \ref{tab:ErrorRatio_w} reveals several insights. The LAH performs best with average error ratios very close to zero, and its performance improves further as capacity in the service region increases. On the other hand, MH performs consistently poorly, with average error ratios varying from 73.6\% to 78.3\%. We see in Figure \ref{fig:Heuristics_w}b that the DP and the LAH ship much more food than the MH. Since the MH does not consider future periods, it ships only the minimal amount in the current period to ensure that overshipment will not occur. However, in the long run, this results in high holding costs. The LAH can anticipate this holding cost and ends up shipping the food that is not necessary to store for the next period. The DP objective function value decreases as $w$ increases, as it is easier to achieve equity in a high-capacity system than a low-capacity one. The LAH also performs better as capacity increases as it is less likely to overship food. The opposite is true for the MH as in a low-capacity system, as its conservative nature reduces overshipments. We should note that in some instances the LAH performs better than the DP. This happens because although we use the DP as a benchmark, it is in fact also an approximation due to the need to discretize the state and action spaces, which may result in suboptimal solutions. 
The BH results in error ratios around 30\% and ships an amount of food between the MH and the LAH as expected. 

\begin{figure}[h!]
\centering
\begin{tabular}{ll}
\subfloat[]{\includegraphics[width=3.2in]{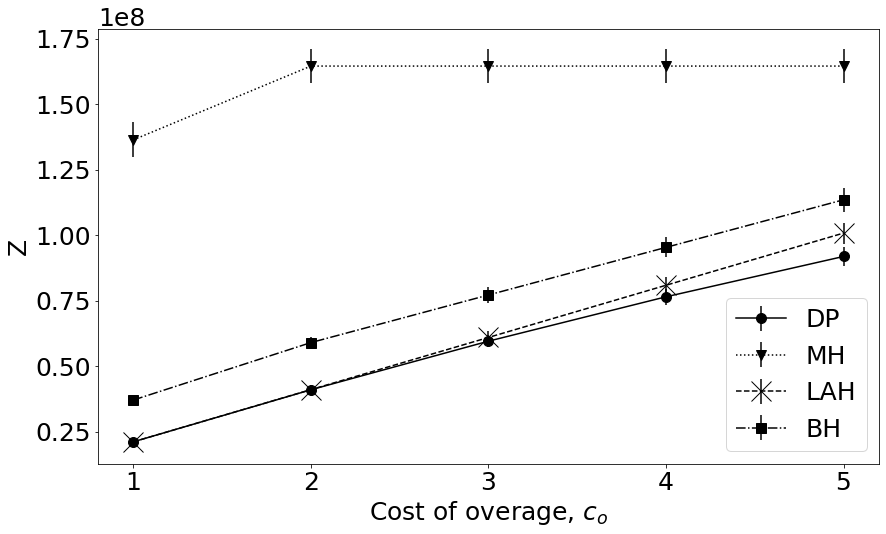}} &
\subfloat[]{\includegraphics[width=3.2in]{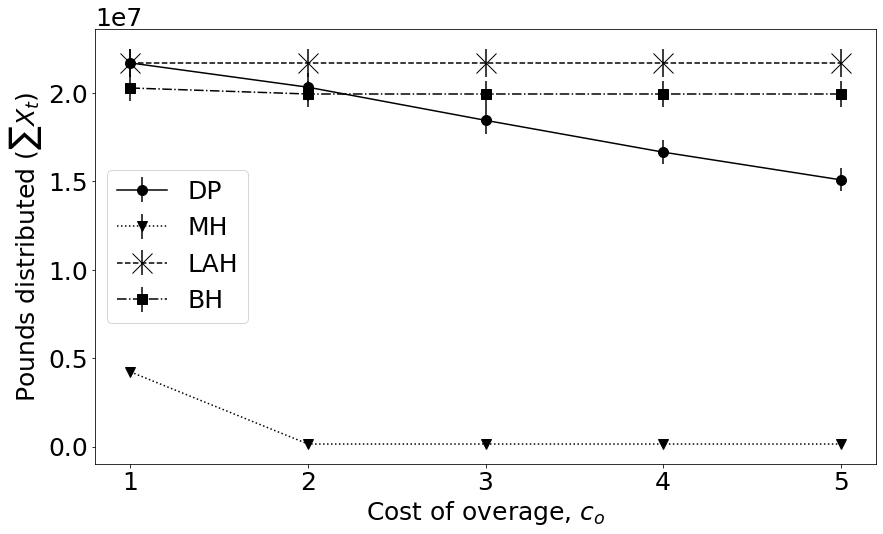}}
\end{tabular}
\caption{Results for (a) the objective function value and (b) the total food distribution over the time horizon for varying $c_o$.}
\label{fig:Heuristics_co}
\end{figure}

\noindent \emph{\textbf{Cost of Overage, $c_o$.}} Examining Figure \ref{fig:Heuristics_co}a together with Table \ref{tab:ErrorRatio_co} shows that the LAH again performs close to optimal but its performance gets worse as $c_o$ increases. This is because as $c_o$ increases, there is higher risk associated with shipping too much, which the DP mitigates by shipping less (see Figure \ref{fig:Heuristics_co}b). However, the LAH cannot anticipate this high cost and ships too much. In contrast, the MH performs better as $c_o$ increases with error ratios decreasing from an average of 43\% when $c_o=1$ to an average of 19\% when $c_o=5$. This is because the MH's conservative nature works well under a high cost of overshipping.

\begin{table}[htbp]
    \centering
    \small
    \setlength{\tabcolsep}{3pt}

    \caption{Statistics (mean (Avg), standard deviation (StD),
    minimum (Min), and maximum (Max)) for the error ratio
    for varying $c_o$.}
    \label{tab:ErrorRatio_co}

    \begin{adjustbox}{max width=\linewidth}
    \begin{tabular}{l|ccc|ccc|ccc|ccc|ccc}
    & \multicolumn{3}{c|}{$c_o=1$}
    & \multicolumn{3}{c|}{$c_o=2$}
    & \multicolumn{3}{c|}{$c_o=3$}
    & \multicolumn{3}{c|}{$c_o=4$}
    & \multicolumn{3}{c}{$c_o=5$} \\
    \hline
    & MH & LAH & BH
    & MH & LAH & BH
    & MH & LAH & BH
    & MH & LAH & BH
    & MH & LAH & BH \\
    \hline
    Avg &84.3\%&0\%&42.9\%
        &75\%&0.1\%&30.6\%
        &63.7\%&2.3\%&22.9\%
        &53.4\%&5.4\%&19.8\%
        &44\%&8.8\%&19.1\%\\
    StD &2\%&0\%&2\%
        &2.4\%&0.5\%&1.7\%
        &3.2\%&2.1\%&1.1\%
        &3.6\%&3.5\%&1.2\%
        &3.9\%&4.5\%&2\%\\
    Min &77.6\%&0\%&37.8\%
        &66.4\%&-1.2\%&25.6\%
        &53.3\%&-0.7\%&20.2\%
        &43\%&-1.2\%&17.4\%
        &33.2\%&0.7\%&15.3\%\\
    Max &87.8\%&0\%&47\%
        &79.6\%&1.3\%&33.7\%
        &69.7\%&8.5\%&25.3\%
        &59.2\%&14.5\%&23.6\%
        &52.1\%&18.1\%&23.7\%\\
    \hline
    \end{tabular}
    \end{adjustbox}
\end{table}

\begin{figure}[h!]
\centering
\begin{tabular}{ll}
\subfloat[]{\includegraphics[width=3.2in]{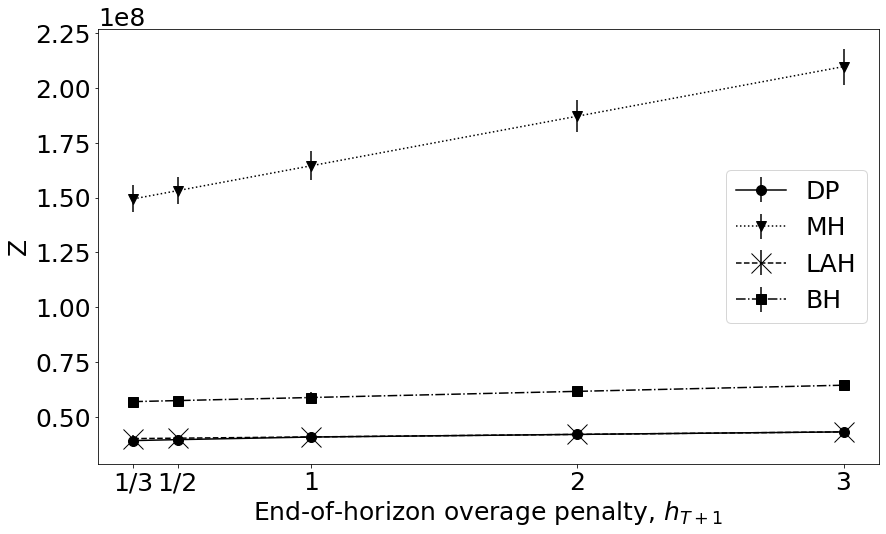}} &
\subfloat[]{\includegraphics[width=3.2in]{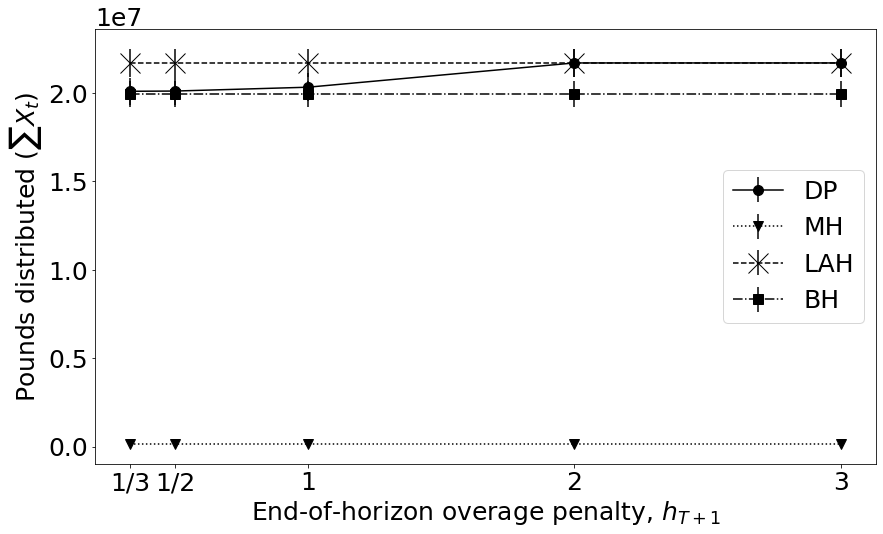}}
\end{tabular}
\caption{Results for (a) the objective function value and (b) the total food distribution over the time horizon for varying $h_{T+1}$.}
\label{fig:Heuristics_hterm}
\end{figure}

\noindent \emph{\textbf{End-of-Horizon Overage Penalty, $h_{T+1}$.}} We obtain similar results when we vary the end-of-horizon overage penalty value. As seen in Figure \ref{fig:Heuristics_hterm} and Table \ref{tab:ErrorRatio_hterm}, the LAH again performs well, improving as $h_{T+1}$ increases. This is expected as the LAH tends to ship more and performs well when overage is more costly. In contrast, MH performs even worse as \(h_{T+1}\) increases. BH performs well in terms of the total shipments over the horizon (see Figure \ref{fig:Heuristics_hterm}b), but still results in error ratios around 30\% mainly due to overshipping in the last period as the BH cannot plan for a high end-of-horizon penalty.

\begin{table}[htbp]
    \centering
    \small
    \setlength{\tabcolsep}{3pt} 

    \begin{adjustbox}{max width=\linewidth}
    \begin{tabular}{l|ccc|ccc|ccc|ccc|ccc}
    & \multicolumn{3}{c|}{$h_{T+1}=1/3$}
    & \multicolumn{3}{c|}{$h_{T+1}=1/2$}
    & \multicolumn{3}{c|}{$h_{T+1}=1$}
    & \multicolumn{3}{c|}{$h_{T+1}=2$}
    & \multicolumn{3}{c}{$h_{T+1}=3$} \\
    \hline
    & MH & LAH & BH
    & MH & LAH & BH
    & MH & LAH & BH
    & MH & LAH & BH
    & MH & LAH & BH \\
    \hline
    Avg &73.5\%&2.2\%&31.1\%
        &73.9\%&1.6\%&30.9\%
        &75\%&0.1\%&30.6\%
        &77.4\%&0\%&31.8\%
        &79.3\%&0\%&33.1\%\\
    StD &2.6\%&2.1\%&1.2\%
        &2.5\%&1.6\%&1.3\%
        &2.4\%&0.5\%&1.7\%
        &2.1\%&0\%&1.4\%
        &2.1\%&0\%&1.7\%\\
    Min &65.4\%&-0.9\%&27.8\%
        &65.6\%&-0.9\%&27.1\%
        &66.4\%&-1.2\%&25.6\%
        &69.8\%&0\%&27.9\%
        &72.2\%&0\%&29.4\%\\
    Max &78.2\%&8.5\%&33.6\%
        &78.6\%&6.4\%&33.7\%
        &79.6\%&1.3\%&33.7\%
        &81.7\%&0\%&34.7\%
        &83.3\%&0\%&36.6\%\\
    \hline
    \end{tabular}
    \end{adjustbox}

    \caption{Statistics (mean (Avg), standard deviation (StD),
    minimum (Min), and maximum (Max)) for the error ratio
    for varying $h_{T+1}$.}
    \label{tab:ErrorRatio_hterm}
\end{table}

\noindent \emph{\textbf{Target Maximum Ending Inventory Level, $\Theta$.}} We find the results are relatively insensitive to changes in $\Theta$, and the heuristics' performance remains consistent with our earlier findings. We therefore present the detailed analyses for varying $\Theta$ in Section~\ref{Appendix:VaryingTheta} of the supplementary material.

\section{Conclusions}\label{sec:conclusions}
%
In this paper, we study the problem of the equitable and effective distribution of food donations from the food bank to the counties in their service region, considering stochastic capacity and supply over a multi-period planning horizon. First, we study the single-period version of this problem and show its equivalence to the newsvendor problem. We use these results to develop a Myopic Heuristic that optimizes shipment decisions in individual periods without considering later periods. We show that this heuristic provides a lower bound on optimal shipments. We then explore the structural properties of the problem and formulate a DP to obtain exact solutions. We use these structural results to obtain a One-Period Look-Ahead Heuristic and a Bounding Heuristic. We perform an extensive numerical study that is calibrated using real-world data from our food bank partner to explore the performance of the heuristics. Our results show that the Look-Ahead Heuristic provides the best performance with an average relative error of 1.1\% among all our experiments. The Myopic Heuristic performs poorly in most cases because it considers only one period at a time. However, its performance improves as the overshipment cost increases, since it tends to ship less than the optimal solution. The Bounding Heuristic generally achieves an intermediate level of performance relative to the other two heuristics.

Our analyses reveal several insights that would be useful to managers of humanitarian organizations. First, we find that planning for multiple periods is better than considering one period at a time. However, the marginal benefit of increasing the planning horizon decreases after a point. Extending the planning horizon from one to eight periods substantially reduces costs, with only modest additional savings beyond eight periods. Second, the overall cost decreases as the capacities of the bottleneck counties increase. So, food banks should focus on increasing the minimum CD ratios in their region to improve their overall performance. 

Future research directions include considering multiple food groups and their interactions in meeting beneficiaries’ nutritional needs while maintaining equitable distribution. The heuristics proposed in this paper could also be enhanced through a learning system where the historical data is collected at the end of each period and used to estimate capacity and supply probability distributions in future periods. We would also like to understand how the results in this paper 
could be extended to consider food perishability constraints. Although the models in this paper have been motivated by a problem that our partner food bank was facing, the results and the developed heuristics can be used to address other problems where equitable and effective distribution of scarce resources are considered in a stochastic capacity and supply environment. Some examples include optimal vaccine allocation before flu season, allocation of teachers to school districts, and distribution of aid kits to the affected population for long-term disasters such as war and drought. 

\section*{Acknowledgments}
This research has been supported by the National Science Foundation under grant numbers 1000018, 1000828, and 2125600. 
The opinions expressed in the paper represent those of the authors and not necessarily those of the National Science Foundation. We would also like to thank the Food Bank of Central and Eastern North Carolina for their support.

\section*{CRediT authorship contribution statement}
\textbf{Irem Sengul Orgut}: Conceptualization, Data curation, Formal analysis, Funding acquisition, Investigation, Methodology, Project administration, Resources, Software, Validation, Visualization, Writing – original draft, Writing – review \& editing. \textbf{Julie Ivy, Reha Uzsoy}: Conceptualization, Funding acquisition, Investigation, Methodology, Project administration, Supervision, Resources, Writing – review \& editing. \textbf{James Wilson}: Investigation, Methodology, Project administration, Resources, Supervision, Writing – review \& editing.

\section*{Appendices}
\appendix

\numberwithin{equation}{section}
\numberwithin{figure}{section}
\numberwithin{table}{section}

\titleformat{\section}
  {\normalfont\large\bfseries}
  {Appendix \thesection.}
  {0.5em}
  {}
\section{Proof of Lemma \ref{lem:recourse}.} \label{Appendix:1}
\begin{proof}



We begin by determining the optimal overshipment quantities $Y_j$ after all shipments $X_j$ have occurred and each county's capacity $\hat{C_j}$ has been observed. Equation \eqref{con:equitySPM2} can be rewritten as
\begin{equation}\label{con:equitySPM2v2}
    Y_j = X_j - \frac{D_{j}(\tilde{X}-\tilde{Y})}{\Delta} \quad \quad \forall j \in \mathcal{J}  \\
\end{equation}

\noindent Inserting (\ref{con:equitySPM2v2}) into (\ref{con:capSPM}), we obtain
\begin{equation}\label{con:capSPM2v2}
   X_j - \frac{D_{j}(\tilde{X}-\tilde{Y})}{\Delta} \geq X_j - C_j \quad \quad \forall j \in \mathcal{J}  \\
\end{equation}
which leads to
\begin{equation}\label{con:capSPM2v3}
   \tilde{Y} \geq \tilde{X} - \Delta \min_j\left(\frac{C_j}{D_j}\right) = \tilde{X} - \hat{R}\Delta\\
\end{equation}
by the definition of the minimum CD ratio, $R$.
Similarly, since by definition $Y_j \geq 0$ for all $j \in \mathcal{J}$, by using (\ref{con:equitySPM2v2}) we obtain
\begin{equation}\label{con:nonnegSPMv2}
   X_j - \frac{D_{j}(\tilde{X}-\tilde{Y})}{\Delta} \geq 0 \quad \quad \forall j \in \mathcal{J}  \\
\end{equation}
which leads to
\begin{equation}\label{con:nonnegSPMv3}
   \tilde{Y} \geq \tilde{X} - \Delta \min_j\left(\frac{X_j}{D_j}\right).\\
\end{equation}

\noindent Combining (\ref{con:capSPM2v3}) and (\ref{con:nonnegSPMv3}), since the objective aims to minimize $\tilde{Y}$, we get
\begin{equation}\label{con:YtildeSS}
   \tilde{Y^*} = \tilde{X} - \Delta \min\left\{\hat{R}, \min_j\left(\frac{X_j}{D_j}\right)\right\},\\
\end{equation}
which, by (\ref{con:equitySPM}) is equal to
\begin{equation}\label{con:YtildeSS}
   \tilde{Y^*} = \tilde{X} - \Delta \min\left\{\hat{R}, \frac{\tilde{X}}{\Delta}\right\}=\tilde{X} - \min\{\hat{R}\Delta, \tilde{X}\} =\max\{0,\tilde{X}-\hat{R}\Delta\}.
\end{equation}
Finally, inserting this result into (\ref{con:equitySPM2v2}), we obtain
\begin{align}
    Y_{j} = X_{j} - \frac{D_{j}(\tilde{X}-\max\{0,\tilde{X}-\hat{R}\Delta\})}{\Delta} = X_{j} - \frac{D_{j} \min\{\tilde{X}, \hat{R}\Delta\}} {\Delta} \\ = X_{j} - \min\{\tilde{X} \frac{D_{j}}{\Delta} , \hat{R}D_j\} = \max\{0, X_j - \hat{R} D_j\}.
\end{align}
The proof is completed.
\end{proof}

\section{Proof of Proposition \ref{prop:MHLB}.} \label{Appendix:3}
\begin{proof}

We start by proving that the functions $U_t(\beta_t,\tilde{X}_t)$ and $V_t(\beta_t,\tilde{X}_t)$ are strictly convex in $\tilde{X}_t$ on the interval $[0,\beta_t]$. The strict convexity of $U_t(\beta_t,\tilde{X}_t)$ has been proven as a part of the proof of Proposition \ref{prop:SPMX}. To prove the strict convexity of $V_t(\beta_t,\tilde{X}_t)$, we will first prove the result for $t=T$, then show that it holds for all $1 \leq t \leq T$. 

First, note that for any $t$ and $\tau$ such that $1 \leq t < \tau \leq T+1$, we have
\begin{equation} \label{eq:Betat}
    \beta_{\tau} = \beta_{t} - \sum_{i=t}^{\tau-1}{\tilde{X}_i} + \sum_{i=t+1}^\tau{S_i}
\end{equation}
so,
\begin{equation} \label{eq:BetaT1D1}
    \frac{\partial \beta_{\tau}}{\partial \tilde{X}_t} = -1.
\end{equation}

\noindent \textbf{Strict convexity of $V_T(\beta_T,\tilde{X}_T)$.} 
Let us write the terminal cost function as:
\begin{equation}
    V_{T+1}(\beta_{T+1}) = \left\{ 
        \begin{array}{ c l }
        h_{T+1}(\beta_T-\tilde{X}_T+S_{T+1}-\Theta) & \textrm{if } S_{T+1} > \Theta + \tilde{X}_T - \beta_T \\
        0 \quad \textrm{otherwise }
        \end{array}
    \right.
\end{equation}

Then, we have
\begin{multline}\label{eq:VT}
    V_{T}(\beta_{T},\tilde{X}_T) = h(\beta_T-\tilde{X}_T) + c_o \tilde{X}_T - c_o \int_{r_L}^{\frac{\tilde{X}_T}{\Delta}} r \Delta f_R(r) dr  - c_o \int_{\frac{\tilde{X}_T}{\Delta}}^{r_U} \tilde{X}_T f_R(r) dr \\ + \int_{\Theta + \tilde{X}_T - \beta_T}^{\infty} h_{T+1} (\beta_T-\tilde{X}_T+s_{T+1}-\Theta) f_{S_{T+1}}(s_{T+1}) ds_{T+1}
\end{multline}

Taking the derivative of $V_{T}(\beta_{T},\tilde{X}_T)$ with respect to $\tilde{X}_T$, we get
\begin{equation}\label{eq:VTDer1}
    \frac{\partial V_{T}(\beta_{T},\tilde{X}_T)}{\partial \tilde{X}_T} = -h + c_o - c_o \left(1-F_R\left(\frac{\tilde{X}_T}{\Delta}\right)\right) - h_{T+1} (1- F_{S_{T+1}}(\Theta + \tilde{X}_T - \beta_T)).
\end{equation}

Now, taking the second derivative, we get:
\begin{equation} \label{eq:VTDer1}
    \frac{\partial^2 V_{T}(\beta_{T},\tilde{X}_T)}{\partial \tilde{X}_T^2} = \frac{c_o}{\Delta} f_R(\frac{\tilde{X}_T}{\Delta}) + h_{T+1} f_{S_{T+1}}(\Theta + \tilde{X}_T - \beta_T)).
\end{equation}

Since by definition, $f_R(.)>0$ and $f_S(.)>0$, we get $\frac{\partial^2 V_{T}(\beta_{T},\tilde{X}_T)}{\partial \tilde{X}_T^2} >0$, proving that $V_{T}(\beta_{T},\tilde{X}_T)$ is strictly convex. 

\noindent \textbf{Strict convexity of $V_t(\beta_t,\tilde{X}_t)$.} Now, let us consider a general $t$ such that $1 \leq t \leq T$ and write $V_{t}(\beta_{t},\tilde{X}_t)$ in its extended form. For the sake of simplicity, we will use the following notation: $E_R[c_o \tilde{Y_t}] = c_o \tilde{X}_t - c_o \int_{r_L}^{\frac{\tilde{X}_t}{\Delta}} r \Delta f_R(r) dr - c_o \int_{\frac{\tilde{X}_t}{\Delta}}^{r_U} \tilde{X}_t f_R(r) dr$. Then,
\begin{multline}\label{eq:VtExtend}
    V_{t}(\beta_{t},\tilde{X}_t) = h(\beta_t-\tilde{X}_t) + E_R[c_o \tilde{Y}_t] + 
    \int_{0}^{\infty} (h(\beta_{t+1}-\tilde{X}_{t+1}) + E_R[c_o \tilde{Y}_{t+1}] \\ + \int_{0}^{\infty} (h(\beta_{t+2}-\tilde{X}_{t+2}) + E_R[c_o \tilde{Y}_{t+2}] + \dots + \int_{0}^{\infty} (h(\beta_{T}-\tilde{X}_{T}) + E_R[c_o \tilde{Y}_{T}]\\
    + \int_{\Theta + \tilde{X}_T - \beta_T}^{\infty} (h_{T+1} (\beta_T-\tilde{X}_T+s_{T+1}-\Theta) f_{S_{T+1}}(s_{T+1}) ds_{T+1} )  f_{S_{T}}(s_{T}) ds_{T}) \\
    \dots f_{S_{t+2}}(s_{t+2}) ds_{t+2} )) f_{S_{t+1}}(s_{t+1}) ds_{t+1})
\end{multline}

Note that $\frac{\partial E_R[c_o \tilde{Y_{\tau}}]}{\partial \tilde{X}_t} = 0$ for $1 \leq t < \tau \leq T$. So, we can write (\ref{eq:VtExtend}) as:
\begin{align}\label{eq:Vt2}
\begin{split}
    V_{t}(\beta_{t},\tilde{X}_t) = h(\beta_t-\tilde{X}_t) + E_R[c_o \tilde{Y_t}] +\\ 
    \int_{0}^{\infty} (h(\beta_{t+1}-\tilde{X}_{t+1}) + E_R[c_o \tilde{Y_{t+1}}]) f_{S_{t+1}}(s_{t+1}) ds_{t+1}+\\
    \int_{0}^{\infty} \int_{0}^{\infty} (h(\beta_{t+2}-\tilde{X}_{t+2}) + E_R[c_o \tilde{Y_{t+2}}]) f_{S_{t+2}}(s_{t+2}) ds_{t+2} f_{S_{t+1}}(s_{t+1}) ds_{t+1}+\\ 
    \dots \\
    \int_{0}^{\infty} \int_{0}^{\infty} \dots \int_{0}^{\infty} (h(\beta_{T}-\tilde{X}_{T}) + E_R[c_o \tilde{Y_{T}}]) f_{S_{T}}(s_{T}) ds_{T} \\ \dots  f_{S_{t+2}}(s_{t+2}) ds_{t+2} f_{S_{t+1}}(s_{t+1}) ds_{t+1} +  \\
    + \int_{0}^{\infty} \int_{0}^{\infty} \dots \int_{0}^{\infty} \int_{\Theta + \tilde{X}_T - \beta_T}^{\infty} h_{T+1} (\beta_T-\tilde{X}_T+s_{T+1}-\Theta) f_{S_{T+1}}(s_{T+1}) ds_{T+1} f_{S_{T}}(s_{T}) ds_{T} \\\dots f_{S_{t+2}}(s_{t+2}) ds_{t+2} f_{S_{t+1}}(s_{t+1}) ds_{t+1}.
\end{split}
\end{align}

\noindent Let us denote the last term of (\ref{eq:Vt2}) as $\Upsilon$ and take its derivative with respect to $\tilde{X}_t$:
\begin{multline}\label{eq:UpsilonDer}
    \frac{\partial \Upsilon}{\partial \tilde{X}_t} =
    -h_{T+1} \int_{0}^{\infty} \dots \int_{0}^{\infty} \int_{\Theta +\sum_{i=t}^{T}{\tilde{X}_i} - \sum_{i=t}^T{S_i} - \beta_t}^{\infty} f_{S_{T+1}}(s_{T+1}) f_{S_{T}}(s_{T}) \\
    \dots f_{S_{t+2}}(s_{t+2}) f_{S_{t+1}}(s_{t+1}) ds_{T+1}\ ds _T\dots ds_{t+1}
\end{multline}

Now, taking the derivative of the entire $V_t(\beta_t,\tilde{X}_t)$ expression with respect to $\tilde{X}_t$ using (\ref{eq:UpsilonDer}):

\begin{multline}\label{eq:XtDer2}
    \frac{\partial V_t(\beta_t,\tilde{X}_t)}{\partial \tilde{X}_t} = -h +c_o - c_o \left[1-F_R\left(\frac{\tilde{X}_t}{\Delta}\right)\right] \\
    + \int_{0}^{\infty} -h f_{S_{t+1}}(s_{t+1}) ds_{t+1} \\
    + \int_{0}^{\infty} \int_{0}^{\infty} -h f_{S_{t+1}}(s_{t+1}) f_{S_{t+2}}(s_{t+2}) ds_{t+2} ds_{t+1} \\
    +\dots+
    \int_{0}^{\infty} \int_{0}^{\infty} \dots \int_{0}^{\infty} -h f_{S_{t+1}}(s_{t+1}) f_{S_{t+2}}(s_{t+2}) \dots f_{S_{T}}(s_{T}) ds_{T} \dots ds_{t+2} ds_{t+1}\\
    -h_{T+1} \int_{0}^{\infty} \dots \int_{0}^{\infty} \int_{\Theta +\sum_{i=t}^{T}{\tilde{X}_i} - \sum_{i=t}^T{S_i} - \beta_t}^{\infty} f_{S_{T+1}}(s_{T+1}) f_{S_{T}}(s_{T}) \\
    \dots f_{S_{t+2}}(s_{t+2}) f_{S_{t+1}}(s_{t+1}) ds_{T+1}\ ds _T\dots ds_{t+1}
\end{multline}

Finally, taking the second derivative of $V_t(\beta_t,\tilde{X}_t)$ with respect to $\tilde{X}_t$:
\begin{equation}
     \frac{\partial^2 V_t(\beta_t,\tilde{X}_t)}{\partial \tilde{X}_t^2} = \frac{c_o}{\Delta} f_R\left(\frac{\tilde{X}_t}{\Delta}\right) + h_{T+1} \int_{0}^{\infty} \dots \int_{0}^{\infty} f_{S_{T+1}}(\Theta +\sum_{i=t}^{T}{\tilde{X}_i} - \sum_{i=t}^T{S_i} - \beta_t) f_{S_{T}}(s) \dots f_{S_{t+1}}(s) ds_T\dots ds_{t+1}
\end{equation}

Since by definition, $f_R(.)>0$ and $f_S(.)>0$, we get $\frac{\partial^2 V_{t}(\beta_{t},\tilde{X}_t)}{\partial \tilde{X}_t^2} >0$, proving that $V_{t}(\beta_{t},\tilde{X}_t)$ is strictly convex.

\noindent \textbf{Proving that $\tilde{X}_t^{MH} < \tilde{X}_t^*$ for any $1 \leq t \leq T$.} So far, we have shown that both $U_t(\beta_t,\tilde{X}_t)$ and $V_t(\beta_t,\tilde{X}_t)$ are strictly convex (up) in $\tilde{X}_t$ on the interval $[0,\beta_t]$. Next, we will use this result to provide the proof for Proposition \ref{prop:MHLB}. We start by noting that the myopic solution $\tilde{X}_t^{MH}$ is always feasible for the multi-period problem for period $t$ since by definition $\tilde{X}_t^{MH} \leq \beta_t$. Since both $U_t(\beta_t,\tilde{X}_t)$ and $V_t(\beta_t,\tilde{X}_t)$ are strictly convex in $\tilde{X}_t$ on the interval $[0,\beta_t]$, to prove Proposition \ref{prop:MHLB}, it suffices to show that $\left. \frac{\partial V_t(\beta_t,\tilde{X}_t)}{\partial \tilde{X}_t}\right|_{\tilde{X}_t=\tilde{X}_t^{MH}} < 0$. 
The quantity of interest is
\begin{multline}
    \frac{\partial V_t(\beta_t,\tilde{X}_t)}{\partial \tilde{X}_t} =-h +c_o - c_o \left[1-F_R\left(\frac{\tilde{X}_t}{\Delta}\right)\right] \\
    + \int_{0}^{\infty} -h f_{S_{t+1}}(s_{t+1}) ds_{t+1} \\
    + \int_{0}^{\infty} \int_{0}^{\infty} -h f_{S_{t+1}}(s_{t+1}) f_{S_{t+2}}(s_{t+2}) ds_{t+2} ds_{t+1} \\
    +\dots+
    \int_{0}^{\infty} \int_{0}^{\infty} \dots \int_{0}^{\infty} -h f_{S_{t+1}}(s_{t+1}) f_{S_{t+2}}(s_{t+2}) \dots f_{S_{T}}(s_{T}) ds_{T} \dots ds_{t+2} ds_{t+1}\\
    -h_{T+1} \int_{0}^{\infty} \dots \int_{0}^{\infty} \int_{\Theta +\sum_{i=t}^{T}{\tilde{X}_i} - \sum_{i=t}^T{S_i} - \beta_t}^{\infty} f_{S_{T+1}}(s_{T+1}) f_{S_{T}}(s_{T}) \dots f_{S_{t+2}}(s_{t+2}) f_{S_{t+1}}(s_{t+1}) ds_{T+1}\ ds _T\dots ds_{t+1}
\end{multline}

Inserting $\tilde{X}_t^{MH} = \min\{ \beta_t, \Delta F_R^{-1}(\frac{h}{c_o})\}$ (Proposition \ref{prop:SPMX}) into the above equation:
\begin{multline} \label{eq:VtDeratMH}
    \left. \frac{\partial V_t(\beta_t,\tilde{X}_t)}{\partial \tilde{X}_t}\right|_{\tilde{X}_t=\tilde{X}_t^{MH}} =  \int_{0}^{\infty} -h f_{S_{t+1}}(s_{t+1}) ds_{t+1}  + \int_{0}^{\infty} \int_{0}^{\infty} -h f_{S_{t+1}}(s_{t+1}) f_{S_{t+2}}(s_{t+2}) ds_{t+2} ds_{t+1} \\
    +\dots+
    \int_{0}^{\infty} \int_{0}^{\infty} \dots \int_{0}^{\infty} -h f_{S_{t+1}}(s_{t+1}) f_{S_{t+2}}(s_{t+2}) \dots f_{S_{T}}(s_{T}) ds_{T} \dots ds_{t+2} ds_{t+1}\\
    -h_{T+1} \int_{0}^{\infty} \dots \int_{0}^{\infty} f_{S_{T+1}}\left(\Theta +\sum_{i=t}^{T}{\tilde{X}_i} - \sum_{i=t}^T{S_i} - \beta_t\right) f_{S_{T}}(s_{T}) \dots f_{S_{t+2}}(s_{t+2}) f_{S_{t+1}}(s_{t+1}) ds_{T}\dots ds_{t+2} \ ds_{t+1}\\
    = -h (T-t) -h_{T+1} \int_{0}^{\infty} \dots \int_{0}^{\infty} f_{S_{T+1}}\left(\Theta +\sum_{i=t}^{T}{\tilde{X}_i} - \sum_{i=t}^T{S_i} - \beta_t\right) f_{S_{T}}(s_{T}) \\
    \dots f_{S_{t+2}}(s_{t+2}) f_{S_{t+1}}(s_{t+1}) ds_{T}\dots ds_{t+1}.
\end{multline}
Since, by definition $f_S(.)>0$, this equation is negative when $h>0$. Further, (\ref{eq:VtDeratMH}) is equal to zero when all holding costs are equal to zero, i.e. $h=h_{T+1}=0$. In that case we would have $\tilde{X}_t^* = \tilde{X}_t^{MH}$.
This completes the proof.

\end{proof}


\bibliographystyle{apalike-ejor}
\bibliography{ref.bib}

\clearpage

\titleformat{\section}
  {\normalfont\large\bfseries}
  {\thesection.}
  {0.5em}
  {}

\setcounter{section}{0}
\renewcommand{\thesection}{S\arabic{section}}

\numberwithin{equation}{section}
\numberwithin{figure}{section}
\numberwithin{table}{section}

\setcounter{page}{1}
\renewcommand{\thepage}{S\arabic{page}}

\ifdefined\theHsection
    \renewcommand{\theHsection}{supp.\arabic{section}}
\fi

\section*{Supplementary Material for
``Multi-Period Donation Allocation for Food Banks
under Stochastic Capacity and Supply''}

\section{Histogram of weekly donations.} \label{Appendix:HistDonation}
\begin{figure}[H]
\begin{center}
\includegraphics[scale=0.5]{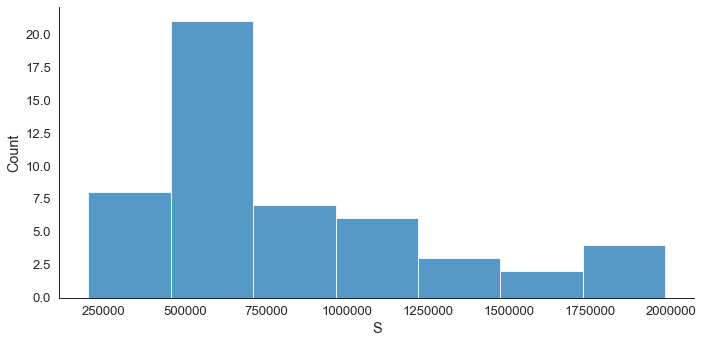}
\caption{Histogram of weekly donations (pounds) to FBCENC during 2018-19 fiscal year.}
\label{fig:SHist}
\end{center}
\end{figure}

\section{Discretization Schemes for Probability Distributions.} \label{Appendix:5}
\subsection{Minimum CD Ratio, $R_t$}

Table \ref{tab:R_Dis} presents the discretization scheme used for $R_t$.

\begin{table}[h!]
    \small
    \centering
    \begin{tabular}{ccc}\hline
    $Q_R$ &	Values &	Probabilities\\ \hline
    2&	0.001604824 & 0.853553\\
	&   0.009353605 & 0.146447\\\hline
    3&	0.00113906 & 0.711093\\
    &	0.006285426 & 0.278518\\
    &	0.017231979 & 0.010389\\\hline
    4&	0.000883655 & 0.603154\\
    &	0.004782699 & 0.357419\\
    &	0.012428557 & 0.038888\\
    &	0.025738805 & 0.000539\\\hline
    5&	0.000722051 & 0.521756\\
    &	0.003872169 & 0.398667\\
    &	0.009852795 & 0.075942\\
    &	0.019412336 & 0.003612\\
    &	0.03463083 & 0.000023\\   
    \end{tabular}
    \caption{Discretization Schemes for Exponential(365.01) Distribution (Miller \& Rice, 1983; Stroud \& Secrest, 1966).}
    \label{tab:R_Dis}
\end{table}

\subsection{Supply, $S_t$}

Table \ref{tab:S_Dis} presents the discretization scheme used for $S_t$.

\begin{table}[h!]
    \small
    \centering
    \begin{tabular}{ccc}\hline
    $Q_S$ &	Values &	Probabilities\\ \hline
    2&	576680.3216&	0.7935709\\
    &	1718743.322&	0.2064291\\\hline
    3&	470542.0216&	0.59771781\\
    &	1249835.822&	0.38307864\\
    &	2728587.022&	0.01920355\\\hline
    4&	409294.2216&	0.458046935\\
    &	1009643.422&	0.469761331\\
    &	2057334.322&	0.070958381\\
    &	3796787.122&	0.001233352\\\hline
    5&	369250.1216&	0.3600101\\
    &	859905.3216&	0.4966017\\
    &	1687398.522&	0.1351428\\
    &	2947960.322&	0.008181715\\
    &	4903192.022&	6.36067E-05\\
    \end{tabular}
    \caption{Discretization Schemes for Gamma Distribution (Miller \& Rice, 1983; Stroud \& Secrest, 1966).}
    \label{tab:S_Dis}
\end{table}

\section{Target Maximum Ending Inventory Level, $\Theta$.}\label{Appendix:VaryingTheta}

\begin{figure}[h!]
\centering
\begin{tabular}{ll}
\subfloat[]{\includegraphics[width=3in]{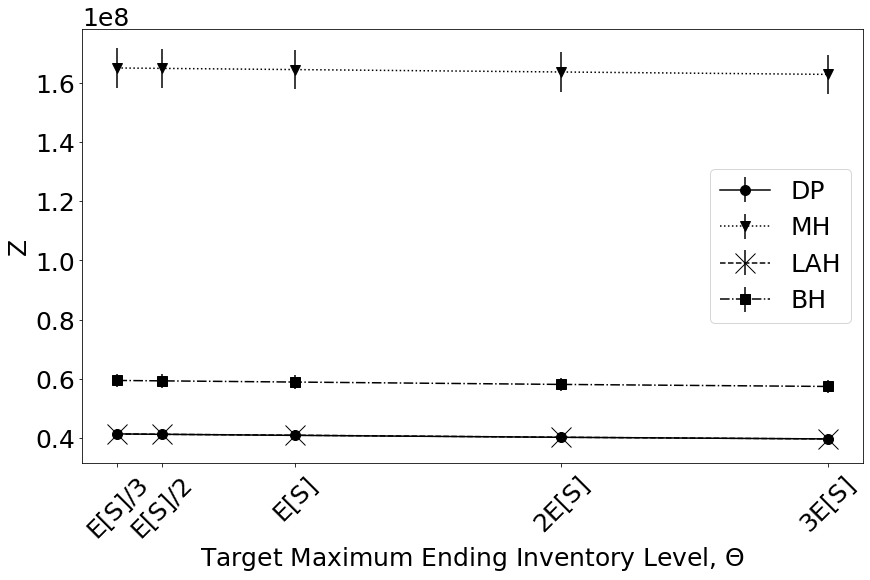}} &
\subfloat[]{\includegraphics[width=3in]{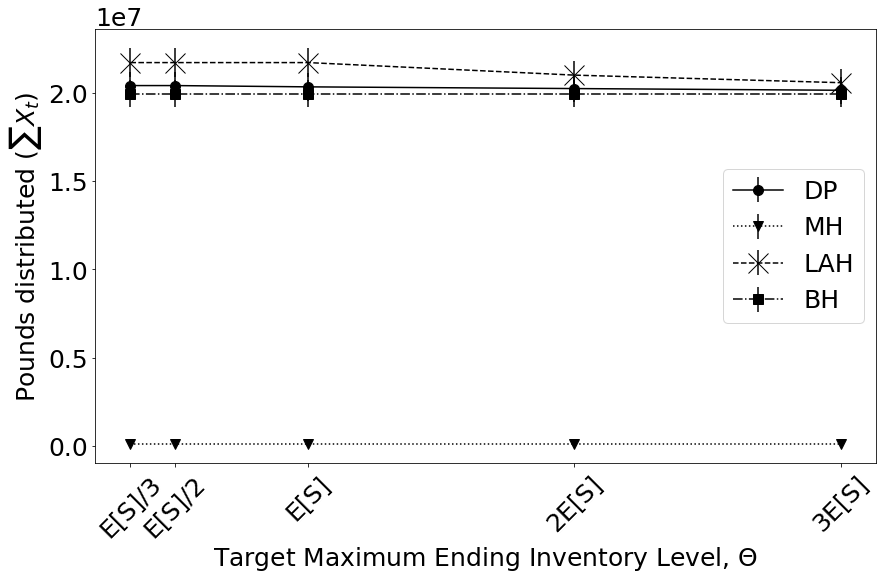}}
\end{tabular}
\caption{Results for (a) the objective function value and (b) the total food distribution over the time horizon for varying $\Theta$.}
\label{fig:Heuristics_Theta}
\end{figure}

We vary the target ending inventory levels and find that results are fairly insensitive to changes in $\Theta$ and are consistent with our findings before. Examining Figure \ref{fig:Heuristics_Theta} and Table \ref{tab:ErrorRatio_Theta} reveals the promising performance of the LAH with the BH resulting in error ratios around 30\% and the MH resulting in error ratios around 75\%.

\begin{table}[htbp]
    \centering
    \small
    \setlength{\tabcolsep}{3pt}

    \caption{Statistics (mean (Avg), standard deviation (StD),
    minimum (Min), and maximum (Max)) for the error ratio
    for varying $\Theta$.}
    \label{tab:ErrorRatio_Theta}

    \begin{adjustbox}{max width=\linewidth}
    \begin{tabular}{l|ccc|ccc|ccc|ccc|ccc}
    & \multicolumn{3}{c|}{$\Theta=E[S]/3$}
    & \multicolumn{3}{c|}{$\Theta=E[S]/2$}
    & \multicolumn{3}{c|}{$\Theta=E[S]$}
    & \multicolumn{3}{c|}{$\Theta=2E[S]$}
    & \multicolumn{3}{c}{$\Theta=3E[S]$} \\
    \hline
    & MH & LAH & BH
    & MH & LAH & BH
    & MH & LAH & BH
    & MH & LAH & BH
    & MH & LAH & BH \\
    \hline
    Avg &74.8\%&0\%&30.4\%
        &74.8\%&0\%&30.4\%
        &75\%&0.1\%&30.6\%
        &75.3\%&0\%&30.7\%
        &75.5\%&0\%&30.9\%\\
    StD &2.4\%&0\%&1.7\%
        &2.4\%&0\%&1.7\%
        &2.4\%&0.5\%&1.7\%
        &2.4\%&0.4\%&1.5\%
        &2.3\%&0.5\%&1.4\%\\
    Min &66.2\%&0\%&25.4\%
        &66.2\%&0\%&25.4\%
        &66.4\%&-1.2\%&25.6\%
        &66.8\%&-1\%&26\%
        &67.3\%&-2\%&26.3\%\\
    Max &79.4\%&0\%&34\%
        &79.5\%&0\%&34.1\%
        &79.6\%&1.3\%&33.7\%
        &79.8\%&1.2\%&33.6\%
        &79.7\%&1.1\%&34.3\%\\
    \hline
    \end{tabular}
    \end{adjustbox}
\end{table}

\end{document}